
\def\LaTeX{%
  \let\Begin\begin
  \let\End\end
  \let\salta\relax
  \let\finqui\relax
  \let\futuro\relax}

\def\UK{\def\our{our}\let\sz s}
\def\USA{\def\our{or}\let\sz z}

\UK



\LaTeX

\USA


\salta

\documentclass[twoside,12pt]{article}
\setlength{\textheight}{24cm}
\setlength{\textwidth}{16cm}
\setlength{\oddsidemargin}{2mm}
\setlength{\evensidemargin}{2mm}
\setlength{\topmargin}{-15mm}
\parskip2mm


\usepackage[usenames,dvipsnames]{color}
\usepackage{amsmath}
\usepackage{amsthm}
\usepackage{mathtools}

\DeclarePairedDelimiter\floor{\lfloor}{\rfloor}
\usepackage{amssymb, bbm}
\usepackage[mathcal]{euscript}
\usepackage[shortlabels]{enumitem}
\usepackage{enumitem}
\usepackage{hyperref}
%
%


\definecolor{viola}{rgb}{0.3,0,0.7}
\definecolor{ciclamino}{rgb}{0.5,0,0.5}
\definecolor{rosso}{rgb}{0.85,0,0}

\def\an #1{{\color{rosso}#1}}
\def\signo #1{{\color{green}#1}}
\def\last #1{{\color{magenta}#1}}

\def\rev #1{{#1}}

\def\last #1{#1}
\def\an #1{#1}
\def\signo #1{#1}
\def\lasts#1{{\color{blue}#1}}



\bibliographystyle{plain}


%
\newtheorem{theorem}{Theorem}[section]
\newtheorem{remark}[theorem]{Remark}

\finqui

\def\Bcenter{\Begin{center}}
\def\Ecenter{\End{center}}
\let\non\nonumber




\def\step #1 \par{\medskip\noindent{\bf #1.}\quad}


\def\Lip{Lip\-schitz}
\def\Holder{H\"older}
\def\Fre{Fr\'echet}

\def\lhs{left-hand side}
\def\rhs{right-hand side}


\def\kk{\kappa}

\def\multibold #1{\def\arg{#1}%
  \ifx\arg\pto \let\next\relax
  \else
  \def\next{\expandafter
    \def\csname #1#1#1\endcsname{{\boldsymbol #1}}%
    \multibold}%
  \fi \next}

\def\pto{.}

\def\multical #1{\def\arg{#1}%
  \ifx\arg\pto \let\next\relax
  \else
  \def\next{\expandafter
    \def\csname #1#1\endcsname{{\cal #1}}%
    \multical}%
  \fi \next}


\def\multimathop #1 {\def\arg{#1}%
  \ifx\arg\pto \let\next\relax
  \else
  \def\next{\expandafter
    \def\csname #1\endcsname{\mathop{\rm #1}\nolimits}%
    \multimathop}%
  \fi \next}

\multibold
qwertyuiopasdfghjklzxcvbnmQWERTYUIOPASDFGHJKLZXCVBNM.

\multical
QWERTYUIOPASDFGHJKLZXCVBNM.


\multimathop
diag dist div dom mean meas sign supp .


\def\Accorpa #1#2 #3 {\gdef #1{\eqref{#2}-\eqref{#3}}%
  \wlog{}\wlog{\string #1 -> #2 - #3}\wlog{}}


\def\supess{\mathop{\rm sup\,ess}}

\def\<#1>{\mathopen\langle #1\mathclose\rangle}
\def\norma #1{\mathopen \| #1\mathclose \|}

\def\[#1]{\mathopen\langle\!\langle #1\mathclose\rangle\!\rangle}

\def\iot {\int_0^t}
\def\ioT {\int_0^T}

\def\intQ{\int_Q}
\def\iO{\int_\Omega}
\def\iG{\int_\Gamma}

\def\dt{\partial_t}
\def\dn{\partial_\nnn}

\def\cpto{\,\cdot\,}

\def\checkmmode #1{\relax\ifmmode\hbox{#1}\else{#1}\fi}

\def\aet{\checkmmode{a.e.\ in~$(0,T)$}}


\def\erre{{\mathbb{R}}}



\def\genspazio #1#2#3#4#5{#1^{#2}(#5,#4;#3)}
\def\spazio #1#2#3{\genspazio {#1}{#2}{#3}T0}

\def\L {\spazio L}
\def\H {\spazio H}
\def\W {\spazio W}

\def\C #1#2{C^{#1}([0,T];#2)}

\def\Hs {\HHH_\sigma}


\def\Lx #1{L^{#1}(\Omega)}
\def\Hx #1{H^{#1}(\Omega)}
\def\LLx #1{\LLL^{#1}(\Omega)}

\def\Wx #1{W^{#1}(\Omega)}

\def\Vs{{\bf V}_\sigma}



\let\th\vartheta
\let\badeps\epsilon
\let\eps\varepsilon
\let\ph\varphi

\let\TeXchi\chi                         
\newbox\chibox
\setbox0 \hbox{\mathsurround0pt $\TeXchi$}
\setbox\chibox \hbox{\raise\dp0 \box 0 }
\def\chi{\copy\chibox}



\let\emb\hookrightarrow

\def\d{\eps}
\def\cd{C_\eps}

\def\vmin{\v_{\rm min}}
\def\vmax{\v_{\rm max}}
\def\J{{\cal J}}
\def\alphaQ{\gamma_1}
\def\alphaO{\gamma_2}
\def\alphav{\gamma_3}
\def\Jred{{\cal J}_{\rm red}}
\def\S{{\cal S}}
\def\Vcon{{\boldsymbol{\cal V}}}
\def\Vad{\Vcon_{\rm ad}}

\def\VR{\Vcon_R}

\def\Vp{V^*}

\let\hat\widehat
\def\ov{\overline}

\let\acc\v
\def\v{{\boldsymbol{v}}}
\def\w{{\boldsymbol{w}}}
\def\con{\v}

\def\vopt{\ov \v}

\def\opt{\vopt}
\def\bph{\ov \ph}
\def\hph{\hat \ph}
\def\bmu{\ov \mu}
\def\hmu{\hat \mu}

\def\optstate{(\bph,\bmu)}
\def\lin{(\xi,\eta)}
\def\adj{(p,q)}
\def\CP{{\bf (CP)}}
\def\Fsec{F''}
\def\ck{{\cal K}}
\def\star{\ast}
\Begin{document}


%
\title{
Regularity results and optimal
velocity control of the 
convective nonlocal Cahn--Hilliard equation in 3D}


%
\date{}
\author{}
\maketitle
\Bcenter
\vskip-1cm
{\large\sc Andrea Poiatti$^{(1)}$}\\
{\normalsize e-mail: {\tt andrea.poiatti@polimi.it}}\\[.25cm]
{\large\sc Andrea Signori$^{(1)}$}\\
{\normalsize e-mail: {\tt andrea.signori@polimi.it}}\\[.5cm]
$^{(1)}$
{\small Dipartimento di Matematica, Politecnico di Milano}\\
{\small via Bonardi 9, 20133 Milano, Italy}



\Ecenter

\Begin{abstract}
\noindent
In this contribution, we study an optimal control problem for the celebrated nonlocal Cahn--Hilliard equation endowed with the singular Flory--Huggins potential in the three-dimensional setting. The control enters the governing {\it state system} in a nonlinear fashion in the form of a prescribed solenoidal, that is a divergence-free, vector field, whereas the {\it cost functional} to be minimized is of tracking-type. The novelties of the {present} paper are twofold: in addition to the control application, the intrinsic difficulties of the optimization problem forced us to first establish new regularity results on the nonlocal Cahn--Hilliard equation that were unknown even without the coupling with a velocity field and are therefore of independent interest. 
\last{This happens to be shown using the recently proved separation property along with \textit{ad hoc} \Holder\ regularities and a bootstrap method.
For the control problem, the existence of an optimal strategy as well as first-order necessary conditions are then established.}

\End{abstract}
\salta
\pagestyle{myheadings}
\newcommand\testopari{\sc Poiatti \ -- \ Signori }
\newcommand\testodispari{{\sc \last{regularity results and optimal
velocity control}}}
\markboth{\testodispari}{\testopari}
\finqui
%

\section{Introduction}
\label{SEC:INTRO}
\setcounter{equation}{0}

Let $\Omega \subset \erre^3$  be some open, bounded domain with smooth
boundary $\Gamma:= \partial \Omega$ and the outward unit normal field $\nnn$. For a prescribed final time $T>0$, we  analyze a suitable optimal control problem, which we are going to present below, for the following initial-boundary value problem:
\begin{alignat}{2}
	\label{eq:gen:1}
	& \dt \ph + \nabla \ph \cdot \v -\div(m(\ph)\nabla \mu)
	=
	0
	\quad && \text{in $Q:=\Omega \times (0,T)$,}
	\\
	\label{eq:gen:2}
	& \mu =  - \badeps\,\ck \ast \ph + \badeps^{-1} F' (\ph) 
	\qquad && \text{in $Q$,}
	\\
	\label{eq:gen:3}
	& (m(\ph)\nabla \mu) \cdot \nnn = 
	0
	\quad && \text{on $\Sigma:=\Gamma \times (0,T)$,} 
	\\
	\label{eq:gen:4}
	& \ph(0)=\ph_0
	\quad && \text{in $\Omega$.}
\end{alignat}
\Accorpa\Sysgen {eq:gen:1} {eq:gen:4}
The system represents a convective version of the celebrated nonlocal Cahn--Hilliard equation, where the velocity field $\v$ is prescribed.
We briefly describe the primary variables of the system as their complete understanding will be clarified below through the presentation of the connected literature.
The variable $\ph$ denotes an order parameter known as the {\it phase field}, $\mu$  is the associated chemical potential, 
$m(\ph)$ is a mobility function, while $\badeps$ is a positive physical constant. Finally, $\v$ stands for a divergence-free vector field that will play the role of control later on, $\ph_0$ is a suitable initial datum, and $F'$ stands for the derivative of a double-well-shaped nonlinearity.


 The above  system represents a nonlocal version of the (local) Cahn--Hilliard  equation that was originally introduced in \cite{num5} and \cite{num13} to model segregation processes in a binary mixture. Despite its original application related to material science, in the last decades the model \last{has proven to be} remarkably flexible \last{in describing} plenty of segregation-driven problems related to cell biology \cite{BH,Dolgin} and tumor growth \cite{CL} (see also \cite{M} and the references therein).
Given  a binary mixture, we indicate with $\ph$ the phase field variable describing the relative mass fraction difference. 
It is assumed that $\{\ph= 1\}:=\{x \in \Omega : \ph(x)= 1\}$ and $\{\ph= -1\}$ indicate the regions occupied by the pure phases and that $\ph$ smoothly transits from $-1$ to $1$ in a narrow transition layer, approximating the interface, whose thickness scales as $\epsilon>0$.
If the mixture is isothermal, the system evolves  minimizing the Ginzburg--Landau functional reading as 
\begin{equation}
    \mathcal{G}(\ph)
    :=\frac{\badeps}{2} \iO |\nabla \ph|^2
    +\frac 1 \badeps \iO \Psi(\ph),
    \label{ener_loc}
\end{equation} 
where $\Psi(\ph)$ is the Flory--Huggins  free energy density
\begin{equation}
    \Psi(s)=\frac{{\theta}}{2}((1+s)\text{ln}(1+s)+(1-s)\text{ln}(1-s))-\frac{\theta_0}{2}s^2=F(s)-\frac{\theta_0}{2}s^2,\ \ \ \forall s\in[-1,1],
    \label{potential_loc}
\end{equation}
with constants related to the mixture temperature such that $0<\theta<\theta_0$. 
As customary, $\Psi$ is also called {\it logarithmic potential} and it is extended by continuity at $\pm1$ \rev{as $+\infty$} otherwise. 
The nonlinearity \eqref{potential_loc} is said to be \textit{singular}, \last{since} it approaches $+\infty$ as its argument tends to the pure phases $\pm1$.
\last{To} simplify the model, one \last{often} considers a  polynomial approximation of  $\Psi$ taking the {\it regular potential} defined as $\Psi_{\rm reg}(s)=\frac{1}{4}(s^2-1)^2$, $s \in \erre$. It is worth recalling that, in the case of polynomial type potentials,  it is not
possible to guarantee the existence of physical solutions, that is, solutions for which $-1 \leq  \ph  \leq  1$ throughout the evolution. Therefore, to stick with the physics of the model, we will concentrate on the singular choice \eqref{potential_loc}.
With the above ingredients, the Cahn--Hilliard equation  reads as follows:
\begin{alignat}{2}
        \label{CHH:1}
        & \partial_t \ph- \div(m(\varphi)\nabla\mu) =0 \quad  &&\text{in } Q,\\
    & \mu=-\badeps\Delta \ph+ \badeps^{-1} \Psi'(\ph)
    \quad &&\text{in } Q,
    \label{CHH:2}
    \\
    & \label{CHH:3}
    \dn \ph = ( m(\ph) \nabla \mu )\cdot  \nnn = 0
    \quad && \text{on $\Sigma$,}
    \\
    & \label{CHH:4}
    \ph(0)=\ph_0 
    \quad && \text{in $\Omega$,}
    \end{alignat}
where the no-flux
condition for the chemical potential $\mu$ entails that no mass flux occurs at the boundary.
Noticing that the free energy $\mathcal{G}$ introduced in \eqref{ener_loc} only focuses on short range interactions between particles,
Giacomin and Lebowitz observed in \cite{giacomin, giacomin1, giacomin2} that a physically more rigorous derivation leads to nonlocal dynamics, formulating the  nonlocal Cahn--Hilliard equation. 
From the modeling, using general approaches of statistical mechanics, the mutual short and long-range interactions between particles are described through convolution integrals weighted by interactions kernels. 
In this case, the gradient term is replaced by a nonlocal spatial interaction integral and the free energy $ \cal G$ is replaced by
\signo{the nonlocal Helmholtz free energy}
\begin{align}
    {\mathcal{E}}(\ph):=-\frac{1}{2}\int_{\Omega\times\Omega} \ck(x-y)\ph(x)\ph(y)\, {\rm dx\, dy}
    +\int_\Omega F(\ph),
    \label{en}
\end{align}
where $\ck$ is a sufficiently smooth symmetric interaction kernel. This functional is characterized by a competition between the mixing entropy $F$, a convex function, and a nonlocal demixing term related to $\ck$. As shown in \cite{giacomin1} (see also \cite{GGG0, GGG, P} and the references therein), the energy $\cal G$ can be seen as an
approximation of $\mathcal{E}$, as long as we redefine $F$ as $\widetilde{F}(x,s)=F(s)-\frac{1}{2}(\ck\ast1)(x)s^2$, $x\in \Omega, s \in[-1,1]$. Indeed, we can rewrite $\mathcal{E}$ as 
\begin{align*}
    \mathcal{E}(\ph)
    & =\frac{1}{4}\int_{\Omega\times \Omega} \ck(x-y)\vert \ph(y)-\ph(x)\vert^2\, {\rm dx\, dy}
    +\int_\Omega F(\ph)
    - \frac 12 \iO {a}\ph^2
    \\ & 
    =\frac{1}{4}\int_{\Omega\times \Omega} \ck(x-y)\vert \ph(y)-\ph(x)\vert^2\, {\rm dx\, dy}
    +\int_\Omega \widetilde{F}(\ph),
\end{align*}
upon setting $a(x):=(\ck\ast 1)(x)$, $x\in\Omega$. We can formally interpret $\widetilde{F}$ as the potential $\Psi$ occurring in  \eqref{ener_loc}, and observe that the (formal) first-order approximation of the nonlocal interaction is $\frac{k}{2}\vert\nabla\ph\vert^2$, for some $k>0$, as long as $\ck$ is sufficiently peaked around zero. In that case, it is also possible to study the 
nonlocal-to-local asymptotics that rigorously justifies the above-sketched intuition: see \cite{Scarpa1,Scarpa2}. 
\rev{In that regard, it is worth noting the recent contributions \cite{CES, ES}, which address a similar issue involving degenerate mobility.}
As we are not focusing on nonlocal-to-local asymptotics, we select the easier formulation with $F$ convex and without \last{the term} $a$, being aware that everything can be straightforwardly reformulated for the other case.
The resulting nonlocal Cahn--Hilliard equation then reads as (see \cite{GGG0,GGG}):
 \begin{alignat}{2}
     \label{nonloc3d:1}
    & \partial_t\ph-\div(m(\ph)\nabla \mu)=0\quad && \text{in }Q,\\
    \label{nonloc3d:2}
    & \mu=-\badeps\,\ck\ast \ph + \badeps^{-1}  F^\prime(\ph)\quad && \text{in }Q,\\
    \label{nonloc3d:3}
    & (m(\ph)\nabla \mu) \cdot \nnn=0\quad && \text{on }\Sigma,\\
    \label{nonloc3d:4}
    & \ph(0)=\ph_0\quad && \text{in }\Omega.
 \end{alignat}
 %
The existence of weak solutions, the uniqueness, and the
existence of the connected global attractor 
  were discussed in \cite{Frigeri,FG, FG1}. For the local case, without the claim of being exhaustive, we refer to \cite{M, Wu} and the references therein. 
\last{In the} the local case, few results are known for degenerate mobilities, i.e., that vanish in the pure phases, \last{namely \signo{just}} the existence of weak solutions obtained in \cite{EG}. The coupling of the nonlocal Cahn--Hilliard equation with logarithmic potential with suitably degenerate mobility has \last{instead} proven very effective. Roughly speaking, the idea is to choose the mobility function in such a way that the two degeneracies, the one of the  singular potential and \last{the one }of the  mobility, compensate  providing \eqref{nonloc3d:1}-\eqref{nonloc3d:2} of a  parabolic structure. This also opens the path to \last{obtain} the existence of strong solutions and continuous dependence results: we refer to \cite{Fmob, Frig, FGGS,  FGK} for more details. However, in contraposition to the local case, less is known for regularity theory in the nonlocal case with constant mobility and logarithmic potential. 
\last{Hereby, we aim at filling this gap by providing new regularity theory for the solutions to \signo{\eqref{eq:gen:1}-\eqref{eq:gen:4}}}
Further results concerning well-posedness and regularity of weak solutions for the nonlocal case are studied in \cite{GGG0}, where the validity of the strict separation property in dimension two for system \eqref{nonloc3d:1}-\eqref{nonloc3d:3} with constant mobility and singular potential were established. This means that if the initial state $\ph_0$ is not a pure
phase, i.e., nor $\ph_0\equiv1$ nor $\ph_0\equiv-1$, then the corresponding solution stays away from the pure states in \last{arbitrarily small positive} time, uniformly with respect to the initial datum. This property in dimension two was crucial to derive further regularity
results as well as the existence of regular finite dimensional attractors, whereas in $3$D only the existence of the (possibly infinite-dimensional) global attractor was proven. In the same work\signo{,} the convergence of a weak solution to a single equilibrium was shown as well. Then, in \cite{GGG}, the same authors propose an alternative argument to prove the strict separation property in dimension two, relying on De Giorgi's iteration scheme. More recently, a similar approach was successfully adopted in \cite{P} by the first author to prove for the first time the validity of the instantaneous strict separation property in dimension three, under weaker assumptions on the singular potential $F$ (cf. the forthcoming assumption \ref{h3}). In particular, it was shown that it is not necessary to assume the growth condition
\begin{align}
    \last{\Fsec(s)\leq Ce^{C\vert F^\prime(s)\vert^\gamma},\qquad \forall s\in(-1,1),\quad \gamma\in(1,2],}
    \label{exp}
\end{align}
for some constant $C>0$. This assumption, fulfilled, e.g., by the logarithmic potential, was essential in \cite{GGG0} and \cite{GGG} for the application of the Trudinger--Moser inequality. \rev{We also mention that assumption \eqref{exp} has recently been weakened in a similar way in the case of two-dimensional local Cahn--Hilliard equation as well (see \cite{GPf}).} Leaning on the result of the separation property \rev{for the nonlocal Cahn--Hilliard equation}, in \cite{P} the author derives extra instantaneous regularization of weak solutions showing that, under suitable assumptions on the interaction kernel $\ck$, from any positive time $\tau>0$ onward, the solution $\varphi$ becomes \Holder\ continuous in space-time and it belongs to $L^\frac 4 3 (\tau,\infty; H^2(\Omega))$. Moreover, it is also shown that any weak solution to \eqref{nonloc3d:1}-\eqref{nonloc3d:3} converges to a single equilibrium. Finally,
it was proved that, given a sufficiently regular initial datum $\ph_0$ which is already strictly separated, the solution $\ph$ strictly separates for any $t\geq0$. To conclude the literature review, we point out \cite{St}, where, by means of a slightly refined version of the proof proposed in  \cite{P} for the validity of the strict separation property in three dimensions, extra regularity for the associated global attractor is \last{proven}.

Let us now move to the control \signo{application}.
To perform the associated  analysis, two mathematical ingredients are essential concerning the solution $\varphi$ to \Sysgen: 
\begin{itemize}
    \item 
    the validity of the strict separation property for any time $t\geq 0$, crucial to deal with the nonlinearity $F$ and its higher-order derivatives.
    \item 
    Extra regularity properties on $\varphi$, important when dealing with continuous dependence estimates in stronger norms.
    \last{Those will be fundamental to show some differentiability properties of the \signo{associated} solution operator.}
\end{itemize}
The former readily follows by \cite[Corollary 4.5, Remark 4.7]{P}, so that our first aim is to show that, assuming $\varphi_0$ and $\v$ smooth enough, there exists a suitably regular solution $\varphi$ to \Sysgen. In particular, extending the 2D result of \cite{AGGP}, we show the existence and uniqueness of a weak and a strong solution $\ph$, according to the regularity of $\ph_0$.

Then, we  \signo{establish}  {additional regularity results for the order parameter.
Namely, provided the initial data and the velocity field are regular enough, we can guarantee the bound
\begin{align*}
   \norma{\ph}_{\L4 {\Hx2} \cap \L3 {\Wx{1,\infty}}}
    \leq C,
\end{align*}
for some $C>0$ depending on the data of the system. Actually, we \signo{also obtained}  intermediate regularity results by \an{adopting minimal} assumptions on the prescribed velocity field.
Taking inspiration from \cite{Frig}, this is achieved by performing a delicate bootstrap argument coupled with a maximal regularity \signo{result of type $L^2$-$L^p$} for parabolic
evolution equations with coefficients
continuous in time (and \Holder\ continuous in space) shown in \cite{Pruss}. }
%
%

As anticipated, the velocity field $\v$ occurring in   \eqref{eq:1} is now considered as a control variable and it is  allowed  to vary in a suitable set, referred to as the {\it control-box}, given by
\begin{align*}
    \Vad: = \{ \v \in \L\infty{L^\infty(\Omega; \erre^3)} \cap\Vcon: \vmin \leq \v \leq \vmax\},
\end{align*}
with $\vmin$ and $\vmax$ given \signo{bounded}  functions, the inequalities being intended componentwise, and where the control space reads as
\begin{align*}
	   \Vcon & := \Big \{ \L2 {L^2(\Omega; \erre^3)} : \div \v =0\text{ in }\Omega,\quad \v \cdot \nnn =0 \text{ on }\Gamma\Big \}.
\end{align*}
The {\it cost functional} we want to minimize is a quadratic-type cost functional and it is defined as 
\begin{align}
	\label{cost}
	\J(\con; \ph)
	=
	\frac {\alphaQ}2 \ioT\iO |\ph - \ph_Q|^2
	+ \frac {\alphaO}2 \int_\Omega |\ph(T) - \ph_\Omega|^2
	+ \frac {\alphav}2 \ioT\iO |\v|^2,
\end{align}
where $\alphaQ, \alphaO,$ and $\alphav$ are nonnegative constants, not all zero, whereas $\ph_Q, \ph_\Omega$ denote some prescribed target functions defined in $Q$ and $\Omega$, respectively.
Optimal control theory for Cahn--Hilliard type systems is rather flourishing and we refer to \cite{HW, ZL}
for optimal control problems related to the classical Cahn--Hilliard equation, and to 
\cite{FGS, FRS} for problems connected to the nonlocal convective Cahn--Hilliard equation when coupled with the Navier--Stokes equation for the velocity field.
For nonlocal Cahn--Hilliard type systems with application to biology, we mention \cite{FLS, RSS} and the references therein (see also \cite{SS} \rev{and \cite{F}}).
In all those scenarios, the control variable linearly enters the system as a given source.

Finally, we point out the very related work
\cite{RS}, where an optimal velocity control problem for the nonlocal convective Cahn--Hilliard equation with degenerate mobility has been addressed. Here, we consider the case of constant mobility and logarithmic potential for which regularity theory was unknown before. 
Moreover, we highlight that the control-box $\Vad$ we consider does not require the control function to be weakly differentiable nor in space nor in time\last{,} which is a much more natural assumption for controls (compare our definition of $\Vad$ with \cite[(1.7)]{RS}).
Let us also refer to \cite{CGS1} and \cite{CGS2} for optimal velocity control problems connected to Cahn--Hilliard type systems.

As they will not play any role in the forthcoming  mathematical analysis, we set for convenience $\badeps=1$ and $m\equiv 1$. Thus, the {\it state system} we are going to study reads as 
\begin{alignat}{2}
	\label{eq:1}
	& \dt \ph + \nabla \ph \cdot \v - \Delta \mu
	=
	0
	\quad && \text{in $Q$,}
	\\
	\label{eq:2}
	& \mu =  - \ck \ast \ph + F' (\ph) 
	\qquad && \text{in $Q$,}
	\\
	\label{eq:3}
	& \dn \mu = 
	0
	\quad && \text{on $\Sigma$,} 
	\\
	\label{eq:4}
	& \ph(0)=\ph_0
	\quad && \text{in $\Omega$.}
\end{alignat}
\Accorpa\Sys {eq:1} {eq:4}
Hence, the optimal control problem we are going to address consists in solving the following minimization problem:\begin{align*}
	\CP
	\quad 
	\min_{\con \in \Vad} \J(\con; \ph),
	\quad \text{subject $\ph$ solves \Sys.}
\end{align*}

\last{
An interesting extension of the cost functional
in \eqref{cost} deals with inserting
$L^1$-penalizations to describe sparsity effects. For instance, in \eqref{cost} one may consider an additional term of the form $\gamma_4 \ioT \iO |\v|$ for a suitable nonnegative constant $\gamma_4$.
This extra term would allow proving sparsity effects for the optimal controls, that is, it may be shown that every optimal control vanishes in a \lq\lq large" region of the parabolic cylinder: see, e.g., \cite{CSS4, GLS_OPT, ST, ST2} for optimal control problems related to phase field systems dealing with sparsity.}

The rest of the paper is {structured this way}: in the next section we set the notation and list our main results. Section \ref{SEC:STATE} is then devoted to the analytical results connected to the state system \Sys. In particular, we provide the aforementioned \last{regularity results} as well as continuous dependence estimates with respect to the prescribed velocity flow.
In the last section, we study the optimal control \CP\ by showing the existence of an optimal strategy and then providing first-order necessary conditions for optimality.

\section{Notation and Main Results}
\label{SEC:RESULTS}
\setcounter{equation}{0}

\subsection{Notation}
\label{SUBSEC:NOT}
Before diving into the mathematical results, let us set some conventions and notations.
First, we indicate with $\Omega$ the spatial domain where the evolution takes place and postulate it to be a bounded and smooth domain of~$\erre^3$ (the lower-space dimensions can be treated easily),
with boundary $\Gamma:=\partial\Omega$. 
Then,  $T>0$ denotes a given final time and we set
\begin{align*}
	Q_t:=\Omega\times(0,t),
	\quad 
	\Sigma_t:=\Gamma\times(0,t)
  \,\, \text{ for every $t\in(0,T]$, and}
	\,\,
	Q:=Q_T,
	\quad 
	\Sigma:=\Sigma_T.
\end{align*}
Next, given a Banach space $X$, we use 
$\norma\cpto_X$, $X^*$, and $\<{\cdot},{\cdot}>_X$
to represent, in the order, its norm, its dual space, and the associated duality pairing.
Then,  the classical Lebesgue and Sobolev spaces on $\Omega$, for $1\leq p\leq\infty$ and $k\geq 0$, are indicated by $L^p(\Omega)$ and $W^{k,p}(\Omega)$,
with the standard convention $\Hx k:=\Wx{k,2}$.
The corresponding norms are $\norma{\cdot}_{\Lx p}:= \norma{\cdot}_p$, $\norma{\cdot}_{W^{k,p}(\Omega)}:=\norma{\cdot}_{W^{k,p}}$, and $\norma{\cdot}_{H^{k}(\Omega)}:=\norma{\cdot}_{H^{k}}$. Similar symbols are employed to denote spaces and norms constructed on $Q$, \signo{$\Sigma$ and $\Gamma$}, whereas bold letters are employed for vector-valued spaces.
For instance, for Lebesgue spaces, we have
\begin{align*}
    \LLx p := L^p(\Omega; \erre^3)= \Big\{\fff : \Omega \to \erre^3 : 
    \fff \,\text{measurable}\,,
    \Big(\iO|{\fff}|^p\Big)^{\frac 1p}<  \infty
    \Big\},
    \quad 
    p \geq 1.
\end{align*}
For some particular spaces, let us fix specific shorthands
\begin{align*}
  H & := \Lx2 , 
  \quad  
  V := \Hx1,   
  \quad
  W := \{v\in H^2(\Omega): \ \dn v=0 \,\hbox{ on $\,\Gamma$}\},
  \\
   \HHH  & := \LLx2 , 
   \quad 
   \VVV := \HHH^1(\Omega)= H^1(\Omega;\erre^3),
\end{align*}
as well as the ones for the solenoidal spaces of the velocity field
\begin{align*}
  \textbf{L}_\sigma^p &:= \{\v \in \textbf{L}^p(\Omega) : \div (\v) =0\text{ in }\Omega,\quad \v \cdot \nnn =0 \text{ on }\Gamma\},\quad\text{ for }p\geq2,
  \quad 
   \HHH_\sigma := \LLL^2_\sigma,
  \\
\Vs & := \{\v \in H^1_0(\Omega;\erre^3): \div (\v) =0 \text{ in }\Omega\}.
\end{align*}
We recall that $\textbf{L}_{\sigma }^{p}$ and $\Vs$
correspond to the completion of $\CCC_{0,\sigma }^{\infty }(\Omega)=C_{0,\sigma }^{\infty }(\Omega ;\mathbb{R}%
^{3})$, namely the space of divergence-free vector fields in $\CCC_{0 }^{\infty }(\Omega) = C_{0}^{\infty
}(\Omega ;\mathbb{R}^{3})$, in the norm of $\textbf{L}^{p}(\Omega )$
and $\textbf{V}$, respectively. 
\signo{The above spaces} are endowed by $\norma{\cdot}:= \norma{\cdot}_H= \norma{\cdot}_\HHH, \norma{\cdot}_V, $ $\norma{\cdot}_W$, and $\norma{\cdot}_{\LLL^p},\norma{\cdot}_{\Vs}$, respectively. 
Moreover, we denote the duality product in $V^\ast$ by $\langle\cdot,\cdot\rangle$. We also indicate by $\PPP_\sigma: \HHH\to \HHH_\sigma$ the standard Leray $\textbf{L}^2$-projector onto $\HHH_\sigma$.
In conclusion, we denote by $C^{\alpha}(\ov{\Omega})$, $\alpha\in(0,1)$\signo{,} the spaces of $\alpha$-\Holder\ continuous functions in $\ov{\Omega}$, whereas by $C^{\beta,\gamma}(\ov{Q})$, $\beta,\gamma\in(0,1)$, we refer to the functions which are $\beta$-\Holder\ continuous in space and $\gamma$-\Holder\ continuous in time, respectively. 
Next, for $v \in \Vp$ we set its generalized mean value by
\begin{align*}
    v_\Omega := \frac1 {\vert \Omega\vert}{\<v,1>},
\end{align*}
where the symbol $1$ denotes the constant function in $\Omega$ that assumes the constant value $1$.
This allows us to define $V_{(0)}$ ($H_{(0)}$, respectively) as the space of functions  $\rev{v}\in V$ ($\rev{v}\in H$, respectively) such that ${v}_\Omega=0$ (notice that, being $v$ at least in $H$, $v_\Omega$ is the usual integral mean value), whereas with $V_{(0)}^*$ we denote the space of $\rev{v^*}\in V^*$ such that $\rev{v_\Omega^*}=\langle \rev{v^*},1\rangle=0$.  
Finally, we recall that $H$ will be identified with its dual as usual. Namely, we have the continuous, dense, and compact embeddings:
\begin{align*}
	W \emb V \emb H \emb V^*
\end{align*}
along with the identification 
\begin{align*}
	\<u,v > = \iO uv \quad 
	\forall
	u \in H, v\in V.
\end{align*}
The Laplace operator
\begin{align*}
    {\cal A}_{0}:V_{(0)}\rightarrow V_{(0)}^*
    \quad \text{defined by} \quad
    \left\langle {\cal A}_{0}u,v\right\rangle_{V_{(0)}^*,V_{(0)}} =\iO \nabla u \cdot \nabla v, 
    \quad \text{$v\in V_{(0)}$},
\end{align*}
is a bijective map between $%
V_{(0)}$ and $V_{(0)}^*$. We denote its inverse by $%
\mathcal{N}\signo{:}={\cal A}_{0}^{-1}:V_{(0)}^*\rightarrow V_{(0)}$. As a consequence, for any $v^*\in V_{(0)}^*$%
, we set $\Vert v^*\Vert _{\ast }:=\Vert \nabla \mathcal{N}v^*\Vert
$, which yields a norm in $V_{(0)}^*$, that is
equivalent to the canonical dual norm. In turn, $v \mapsto \big(\Vert v-{v}_\Omega%
\Vert _{\ast }^{2}+|{v}_\Omega|^{2}\big)^{\frac{1}{2}}$ defines a norm in $%
V^*$, that is equivalent to the standard dual norm in $\Vp$. Moreover, by the regularity theory for the
Laplace operator with homogeneous Neumann boundary conditions, there exists a constant $C>0$
such that
\begin{equation}
\Vert \nabla \mathcal{N}v\Vert _{V}\leq C\Vert v\Vert,\quad \forall \,v\in H_{(0)}.  \label{H_2}
\end{equation}%
In conclusion, we introduce the Besov spaces (see, e.g., \cite{Adams} and \cite{Triebel} for more details), as  the following (real) interpolation spaces
$$
    \BB_{p,q}^s(\Omega):=\left(\Lx p,W^{m,p}(\Omega)\right)_{\lambda, q},
$$
where $s=\lambda m$, $m\in\mathbb{N}$, $\lambda\in(0,1)$, $p,q\in[1,\infty)$. In particular, we recall that  
$$
  \BB_{p,q}^{2-\frac 2 q
 }(\Omega)=\left(\Lx p,W^{2,p}(\Omega)\right)_{ 1-\frac 1 q,q},
$$
for any $p,q>1$. 
%
%
\subsection{Assumptions and main results}
Here, we  collect  our main results and the related hypotheses. The following structural assumptions will be in order:
\begin{enumerate}[label={\bf H{\arabic*}}, ref={\bf H{\arabic*}}]
	\item The spatial kernel is such that $\ck \in W_{\rm loc}^{1,1}(\mathbb{R}^3)$, with $\ck({x})=\ck(-{x})$, $x \in \Omega$.
	\label{h1}
	\item The double-well potential fulfills $F\in C^0([-1, 1]) \cap C^3(-1,1)$  and
	\begin{equation*}
	\lim_{s\to-1^+} F'(s)=-\infty, \quad \lim_{s\to1^-} F'(s)=+\infty,\quad F''(s)\geq{\alpha}, \quad\forall\ s\in(-1,1)\rev{, \alpha>0}.
	\end{equation*}
	As usual, we extend $F(s)=+\infty$  for any $s\notin[-1, 1]$.  Without loss of generality, we require $F(0) = 0$ and $F'(0)=0$. In particular, those entail that $F (s) \geq 0$ for any $s \in[-1, 1]$. 
	\label{h2}
 \item \label{h3}
   As $\delta\to 0^+$, we assume
	\begin{align}
	\frac{1}{F^{\prime}(1-2\delta)}={\cal O}\left(\frac{1}{\vert\ln(\delta)\,\vert}\right),\quad \frac{1}{F^{\prime\prime}(1-2\delta)}={\cal O}(\delta),
	\label{F}
	\end{align}
	and analogously that
	\begin{align}
	\dfrac{1}{\vert F^{\prime}(-1+2\delta)\vert }={\cal O}\left(\frac{1}{\vert\ln(\delta)\,\vert}\right),\qquad \dfrac{1}{\Fsec(-1+2\delta)}={\cal O}\left(\delta\right).
	\label{F2}
	\end{align}
	\item \label{h4}
 Either $\ck \in W^{2,1}({B}_\rho)$, where ${B}_\rho:=\{x\in \mathbb{R}^3: \vert x\vert<\rho\}$, with $\rho\sim\text{diam(}\Omega\text{)}$ such that $\overline{\Omega}\subset {B}_\rho$ or $\ck$ is admissible in the sense of \cite[Def.1]{Bedrossian}.
 \label{H4}
\end{enumerate}
\begin{remark}
	As remarked in \cite{Bedrossian}, we observe that Newtonian and Bessel potentials do satisfy assumption \ref{H4}.
\end{remark}

\begin{theorem}
\label{ExistCahn} Let the assumptions \ref{h1}-\ref{h2} \signo{be fulfilled}. Assume that $\v\in L^{4}(0,T;\textbf{L}^6_\sigma)$, 
{
$\ph _{0}\in H$ with $F(\ph_0) \in \Lx1$
}%
and $|(\ph _{0})_\Omega|<1$. Then, there
exists a unique weak solution $(\ph,\mu)$ to \eqref{eq:1}-\eqref{eq:4} in the sense that
\begin{align}
    \label{wreg:1}
\ph & \in \H1 \Vp \cap \C0 H \cap L^{2}(0,T;V), \\
\label{wreg:2}
\ph &\in L^{\infty }(Q):\quad |\ph|<1 \ \text{a.e.
in }Q, \\
\label{wreg:3}
\mu &\in L^{2}(0,T;V),\quad F^{\prime }(\ph )\in
L^{2}(0,T;V),
\end{align}%
and it  satisfies%
\begin{alignat}{2}
\label{weak-nCH-2:1}
    &\< \dt\ph ,v>
    - \iO \ph \v \cdot \nabla v
    + \iO \nabla \mu \cdot \nabla v 
    =0
    \quad && \text{for every $v \in V $, and \aet}, \\
    \label{weak-nCH-2:2}
    & \mu =-\ck\ast \ph + F^{\prime }(\ph ) \quad && \text{a.e. in }Q,%
\end{alignat}
and $\ph ({0})=\ph _{{0}}$ almost everywhere in $\Omega $. The weak
solution fulfills the energy identity
\begin{equation}
    \mathcal{E}(\ph (t))
    +\int_{Q_t}\ph\v\cdot \nabla \mu  
    +\int_{Q_t}|\nabla \mu |^{2} =\mathcal{E}(\ph _{0}),\quad \forall \,t\in \lbrack 0,T].  \label{EE-nCH}
\end{equation}%
In addition, given two weak solutions $\ph_{1}$ and $\ph _{2}$
corresponding to the initial data $\ph _{0}^{1}$ and $\ph _{0}^{2}$ \last{assumed to fulfill the same conditions as above,} and \last{a prescribed} velocity field \last{$\v\in L^{4}(0,T;\textbf{L}^6_\sigma)$}, it holds that
\begin{align}
    \non & 
    \left\Vert \ph_{1} - \ph_{2}\right\Vert _{\C0 \Vp \cap \L2 H} 
    \\ & \quad \non
    \leq \left( \left\Vert \ph _{0}^{1}-\ph _{0}^{2}\right\Vert _{\Vp}+\left\vert (\ph _{0}^{1})_\Omega-(\ph _{0}^{2})_\Omega%
    \right\vert ^{\frac{1}{2}}\Vert \Lambda \Vert _{L^{1}(0,T)}^{\frac{1}{2}%
    }+C \left\vert (\ph _{0}^{1})_\Omega-(\ph _{0}^{2})_\Omega
    \right\vert \right) \times
    \\ & \qquad \times\mathrm{exp} \left(C\left( 1+\Vert \v\Vert
    _{L^{4}(0,T;\textbf{L}^{6}_\sigma)}^{4}\right) \right),
    \label{cd}
\end{align}%
where $\Lambda =2(\left\Vert F^{\prime }(\ph
_{1})\right\Vert _{1}+\left\Vert F^{\prime }\left(\ph _{2}\right)\right\Vert
_{1})$ and $C$ only depend\signo{s} on $\alpha $, $\ck$, $T$, and $\Omega $.
\medskip

\noindent Furthermore, the following regularity results hold:

\begin{itemize}
\item[(i)] \label{i1} If \signo{also} $\ph _{0}\in V$ and it is such
that $F^{\prime }(\ph _{0})\in H$ and $F^{\prime \prime }(\ph
_{0})\nabla \ph _{0}\in \textbf{H}$, then, \last{additionally},
\begin{align}
    \label{reg:2}\ph & \in L^{\infty }(0,T;V)\cap L^{q}(0,T;W^{1,p}(\Omega)),
    \quad q=\frac{4p}{3(p-2)},\quad \forall \,p\in (2,\infty ), \\
    \label{reg:3}
    \dt\ph & \in L^{4}(0,T;\Vp)\cap
    L^{2}(0,T;H), \\
    \label{reg:4}
    \mu & \in L^{\infty }(0,T;V)\cap L^{2}(0,T;W), \\
     \label{reg:5}
    F^{\prime }(\ph ) & \in L^{\infty }(0,T;V).
\end{align}%
and $\partial_\nnn\mu=0$ almost everywhere on $\Sigma$.
 Moreover, if $\v\in L^{\infty }(0,T;\textbf{H}_\sigma)$, we also have $\dt\ph \in L^{\infty
}(0,T;\Vp)$. \smallskip

\item[(ii)] \label{i2} Let the assumptions of $(i)$ hold, together with assumptions \ref{h3} and \ref{h4}. Suppose also that 
\signo{$\ph_0 \in \Lx\infty$ with} $\Vert
\ph _{0}\Vert _{\infty}\leq 1-\delta _{0},$ for some $\delta
_{0}\in (0,1)$. Then, there exists $\delta \signo{\in (0,\delta_0]}$ such that
\begin{equation}
\supess_{t\in \lbrack 0,T]}\Vert \ph (t)\Vert _{\infty}\leq
1-\delta .  \label{SP}
\end{equation}%
As a consequence, we also have that $\mu \in \H1 H \cap \C0 V$.


\item[(iii)] Under the same assumptions of $(ii)$, assume additionally that there exists $\beta_0\in(0,1]$ such that $\varphi_0\in C^{\beta_0}(\overline{\Omega})$. Then there exists $\beta\in(0,\signo{\beta_0]}$ such that
\begin{align}
    \label{holderphi}&\ph\in C^{\beta,\frac \beta 2}(\ov Q),\\&
    \mu\in C^{\beta,\frac \beta 2}(\ov Q),
    \label{holder}
\end{align}
where $\beta$ also depends on the $L^4(0,T;\textbf{L}^6_\sigma)$-norm of $\v$.
\item[(iv)] 
Let the assumptions in point $(iii)$ and  \ref{H4} be fulfilled. Then there exists $C>0$, depending also on the constant $\delta$ appearing in \eqref{SP} and T, such that  
\begin{align}
    \Vert \ph\Vert_{L^2(0,T;H^2(\Omega))}\leq C.
    \label{p2}
\end{align}
Moreover, let us set
\begin{align}
    \label{def:muz}
    \mu_0 := - \ck \ast \ph_0 + F'(\ph_0) .
\end{align}
Then, if $\mu_0\in {\BB}_{3,2}^{1}(\Omega)$, there exists $C>0$ \signo{depending on structural data of the system}, such that 
\begin{align}
    &  \norma{\mu}_{\H1 {\Lx3} \cap \L2 {\Wx{2,3}}}
    + \norma{\ph}_{\H1 {\Lx3} \cap \L4 {\Wx{1,6}}}
    \leq C.
    \label{regul}
\end{align}%
{%
\item[(v)] 
Let the assumptions in point $(iv)$ be fulfilled, $\v \in \L4 {\LLL^\infty}$, and $\mu_0 \in {\cal B}_{6,2}^1\signo{(\Omega)}$. Then, 
there exists $C>0$,  \signo{depending on structural data of the system,} such that 
\signo{\begin{align}
    & \non
    \norma{\mu}_{\H1 {\Lx6} \cap \L3 {\Wx{1,\infty}} \cap\L2 {\Wx{2,6}}} 
    \\ & \quad 
    + \norma{\ph}_{\an{\H1 {\Lx6}\cap}\L3 {\Wx{1,\infty}}}
    \leq C.
    \label{regul:special}
\end{align}}%
\an{Moreover, let instead \signo{assume} $\v\in L^\infty(0,T;\LLL_\sigma^4)$ and $\mu_0\in \BB^{\frac{3}2}_{2,4}(\Omega)$, together with the assumptions in point $(iv)$}.
Then, there exists $C>0$, \signo{depending on structural data of the system,} such that 
\signo{\begin{align}
    & \non
    \norma{\mu}_{\W{1,4} H \cap L^8(0,T;W^{1,4}(\Omega)) \cap \L4 {\Hx2}}
   \\ & \quad 
   +\norma{\ph}_{ L^8(0,T;W^{1,4}(\Omega)) \cap \L4 {\Hx2}}\an{\leq C.}
    \label{reg:v:final}
\end{align}}%
}%
\end{itemize}
\end{theorem}
\begin{remark}
   \last{We \signo{remark} that the assumption $F'(\ph_0)\in L^1(\Omega)$ already implies \signo{that $\ph_0\in\Lx\infty$ with} $\Vert \ph_0\Vert_{\infty}\leq 1$, due to assumption \ref{h2}.} \an{Furthermore, observe that if $\ph_0\in V$ and \signo{it} is strictly separated, i.e., $\Vert
\ph _{0}\Vert _{\infty}\leq 1-\delta _{0},$ for some $\delta
_{0}\in (0,1)$, this directly implies $F'(\ph_0)\in H$ and $F''(\ph_0)\nabla\ph_0\in \HHH$.}
\label{regularity1}
\end{remark}
\begin{remark}
   \last{ We \an{notice} that \cite[Remark 4.4]{AGGP} still holds also in the three-dimensional setting, i.e., any weak solution $\ph$ satisfying \eqref{wreg:1}-\eqref{wreg:2} is instantaneously strictly separated from pure phases.}
    Moreover, under the assumption of the above theorem, \last{part (iii)}, thanks to the \Holder\ regularity in \eqref{holderphi}, the inequality \eqref{SP} reduces to
    \begin{align*}
        \max_{(x,t) \in \ov Q} \norma{\ph(x,t)}_{C^0(\ov Q)}
        \leq 1- \delta.
    \end{align*}
\end{remark}
\begin{remark}
The technical assumption on the initial condition of point $(iv)$ can be avoided by requiring, for instance, that
    $$\mu_0\in W^{2,3}(\Omega)\hookrightarrow \BB_{3,2}^1(\Omega),$$
    which is nevertheless more restrictive. Indeed, from \cite{Adams} we have the embedding
    $$
        \BB_{q,2}^1(\Omega)\hookrightarrow W^{1,q}(\Omega)\hookrightarrow \BB_{q,q}^1(\Omega)\quad \forall q\geq2,
    $$
    and thus the space $\BB_{3,2}^1(\Omega)$ is actually not so far from $W^{1,3}(\Omega)$: Besov spaces set somehow a finer scale with respect to classical Sobolev spaces.
\end{remark}
\begin{remark}
    As it will be clear from the proof, points $(i)$ and $(ii)$ of the above theorem  can be shown 
    by arguing along the same lines of arguments as in \cite{AGGP} which analyzes the same system in the two-dimensional setting.
    The extension to the 3D case we perform  \last{has} been made possible by the recent result of the validity of the strict separation property proven in \cite{P}.
\end{remark}
\begin{remark}  \label{REM:regula}
    Observe that, in the case of point $(iv)$, thanks to the regularity in \eqref{wreg:1}-\eqref{wreg:2}, it holds $\ph\in \C0 V$, the function $t\mapsto \Vert \nabla\ph(t)\Vert^2$ is $AC([0,T])$, and that
\begin{align*}
    \text{$\frac 1 2\frac d {dt}\Vert \nabla\ph(t)\Vert^2=- \iO \dt\ph(t) \Delta\ph(t)$ \quad  for almost every $t\in(0,T)$.}
\end{align*}
\end{remark}
\begin{remark}
\an{Note that, up to point (iv) \signo{it suffices that $\v\in L^4(0,T;\LLL_\sigma^6)$ as long as we assume suitably regular initial data. This} is the minimal summability to get well-posedness of strong solutions and it is enough to deduce \eqref{regul}.}
    \an{Furthermore,}{ using the property $\ph\in L^4(0,T;W^{1,6}(\Omega))$ specified in \eqref{regul}, one may establish further regularity results by \signo{formally} differentiating \eqref{weak-nCH-2:1} in time and testing it by $\partial_t\ph$ as in the 2D analogue (see, e.g., \cite[Thm.2]{FGK} and \cite[Lemma 5.7]{GGG0}). That will readily prove that $\partial_t\ph\in L^\infty(0,T;H)$ and $\mu\in L^\infty(0,T;H^2(\Omega))$. Nevertheless, some extra regularity on the time derivative $\partial_t\v$ \signo{and on the initial data are} required. Being the velocity field \signo{$\v$} our control variable, we do not want to assume $\dt \v \in L^2(Q)$. 
    On the other hand, if one includes that condition in $\Vad$ as done in \cite{RS}, the above regularities can be easily shown.}
\end{remark}
\begin{theorem}\label{THM:CD}
\signo{Suppose that \ref{h1}-\ref{h4} hold.}
\begin{itemize}
    \item[\signo{(i)}]

Let \an{$\v_1$ and $\v_2$} be two given velocity fields such that 
\begin{align*}
    \v_{i} \in \L4 {\LLL^6_\sigma}\signo{,}
    \quad
    {i=1,2}.
\end{align*}
Denote by $(\ph_i,\mu_i)$, $i=1,2$, the two corresponding solutions to \Sys\ related to initial data $\ph_0^i$ which fulfill the assumptions of $(i)$ in Theorem \ref{ExistCahn} and
\rev{$(\ph_{0}^1)_\Omega=(\ph_{0}^2)_\Omega$ with}
\begin{align}
    | {(\ph_{0}^i)_\Omega}| <1,
    \quad 
    \norma{\ph_0^i}_\infty \leq 1- \delta_0^i
    \quad \text{with $\delta_0^i \in (0,1)$,}
    \quad i=1,2,
    \label{conditions}
\end{align}
\signo{exists $\beta_{0,2} \in (0,1]: \ph_0^2 \in C^{\beta_{0,2}}(\ov \Omega)$, and} {$\mu_0^2:= - \ck \ast \ph_0^2 + F'(\ph_0^2)\in \BB_{3,2}^{1}(\Omega)$}. The two solutions $\ph_1$ and $\ph_2$ are then intended in the sense of points $(ii)$ and $(iv)$ of Theorem \ref{ExistCahn}, respectively.
Then, there exists a positive constant $C$ such that 
\begin{align}
    & \norma{\ph_1-\ph_2}_{
	\L\infty H \cap \L2 V}
    \label{cont:dep:est1}
    \leq 
	C 
	(\norma{\v_1-\v_2}_{\L2 \Hs}
 +  \norma{\ph_0^1-\ph_0^2}),
\end{align}
where $C$ depends only on the structure of the system. 

\item[\signo{(ii)}]
{
\noindent \an{Moreover,  suppose that, additionally,
$$
    \v_2\in \L \infty {\LLL^4_\sigma} \cap \L4 {\LLL^\infty},
    \quad \mu_0^1 \in \BB_{3,2}^{1}(\Omega),\quad \mu_0^2 \in \BB_{2,4}^{\frac 32}(\Omega) \cap  \BB_{6,2}^{1}(\Omega),
$$
with \an{$\mu_0^\signo{i}:= - \ck \ast \ph_0^\signo{i}+ F'(\ph_0^\signo{i})$}\signo{, $i=1,2$}\signo{, and exists $\beta_{0,1} \in (0,1]: \ph_0^1 \in C^{\beta_{0,1}}(\ov \Omega)$}.
Then, the two solutions $\ph_1$ and $\ph_2$ are then intended in the sense of points $(iv)$ and $(v)$ of Theorem \ref{ExistCahn}, respectively,} and they fulfill
\begin{align}
    & \norma{\ph_1-\ph_2}_{
	\L\infty V \cap \L2 W}
    \label{cont:dep:est}
    \leq 
	C 
	(\norma{\v_1-\v_2}_{\L6 \Hs}
 +  \norma{\ph_0^1-\ph_0^2}_V)
\end{align}
for a positive constant $C$ which depends only on the structure of the system.
}
\end{itemize}
\end{theorem}
\noindent Once the above analytical properties on the solutions of \Sys\ have been derived, we can address the optimal control problem \CP.  For \an{such a }problem, we postulate the following assumptions:
\begin{enumerate}[label={\bf C{\arabic*}}, ref={\bf C{\arabic*}}]

\item \label{ass:control:1:kernel}
    The spatial kernel $\ck$ fulfills \ref{h1} and \ref{h4}.

\item \label{ass:control:2:initialdata} 
The initial dat{a} fulfill (recall \eqref{def:muz} and Remark \ref{regularity1})
  \an{\begin{align*}
    &
    \ph _{0}\in V \cap \Lx\infty,
    \quad 
    | (\ph_{0})_\Omega| <1,
    \quad 
    \norma{\ph_0}_\infty \leq 1 - \delta_0 \quad \text{with $\delta_0 \in (0,1)$},
    \\ 
    & 
      \signo{\text{exists} \,\, \beta_0 \in (0,1] : \ph_0 \in C^{\beta_0}(\ov \Omega),}
      \quad 
      \mu_0 \in \BB_{2,4}^{\frac 32}(\Omega) \cap  \BB_{6,2}^{1}(\Omega).
\end{align*}}%
\item \label{ass:control:3:const}
    The constants $\alphaQ, \alphaO,$ and $\alphav$ in \eqref{cost} are nonnegative, but not all zero.
    
\item \label{ass:control:4:target}
    The target functions $\ph_Q $ and $\ph_\Omega$ are such that $\ph_Q \in L^2(Q)$,  and $\ph_\Omega \in V$.

\item \label{ass:control:5:Vad}
    The prescribed functions $\vmin$ and $\vmax$ are such that $\vmin, \vmax \in \LLx\infty$ and $\vmin \leq \vmax$ componentwise.

\item \label{ass:control:6:potreg}
    In addition to \ref{h2} and \ref{h3},  the double-well potential is such that $F \in C^4(-1,1)$.
\end{enumerate}
We now present the two main results on \CP. First, we state the existence of an optimal strategy, and then the first-order optimality conditions for minimizers.
\begin{theorem}
    \label{THM:EXCONTROL}
    Assume that \ref{ass:control:1:kernel}-\ref{ass:control:5:Vad} are in force. Then, the optimization problem \CP\ admits at least one solution.
\end{theorem}
\noindent In the formulation of the \signo{corresponding} optimality conditions, we refer to the adjoint variables $p$ an $q$. Those are the unique solutions to a system, referred to as the {\it adjoint system} related to \Sys\ (cf. \eqref{eq:adj:1}-\eqref{eq:adj:4}).
To keep the presentation as essential as possible, we defer to later \last{their} proper introduction (cf. Section \ref{SEC:CONTROL}).
\begin{theorem}
    \label{THM:FOC}
    Assume that \ref{ass:control:1:kernel}-\ref{ass:control:6:potreg} are in force. Let $\opt$ be an optimal control with corresponding state $\optstate$ and adjoint variables $\adj$.
    Then, it necessarily fulfills the variational inequality
    \begin{align}
	\label{foc:final}
	\intQ \big({-}\PPP_\sigma (p \nabla \bph)+ \alphav \vopt\big) \cdot (\v-\vopt)
	\geq 0
	\quad \forall \con \in \Vad.
\end{align}
Moreover, whenever $\alphav \not = 0$, the optimal control $\vopt$ reduces to the $L^2$-orthogonal projection of ${\alphav^{-1} \PPP_\sigma({p \nabla \bph}})$ onto the convex set $\Vad$.
\end{theorem}

\begin{remark}
    Observe that, when $\alphav\not=0$,  \eqref{foc:final} shows that to identify the optimal $\vopt$ is enough to have access, for almost any $t\in(0,T)$, to the divergence-free projection of ${p \nabla \bph}(t)$, i.e., its $\LLL^2$-projection onto $\HHH_\sigma$, whereas its orthogonal complement can be neglected. 
 {We also remark that the standard pointwise characterization of the projection as a suitable bang-bang control involving
 the bounds $\vmin$ and $\vmax$ 
 does not work here, since the result cannot be guaranteed to be divergence-free and thus it may not belong to $\Vad$.}
\end{remark}

Without further reference later on, in the forthcoming \last{estimates} we are going to perform, the capital letter $C$ will denote a generic positive constant that depends only on the structural data of the problem. For this reason, its meaning may change from line to line and even within the same chain of computations.
When it depends on an additional constant $\d$ whose value has to be chosen just at the end of some computations, we use $\cd$  to stress that dependence. 

\section{Mathematical Analysis of the State System}
\label{SEC:STATE}
\setcounter{equation}{0}

\subsection{Proof of Theorem \ref{ExistCahn}}
The proof of the theorem can be mutuated in part from \cite[Thm. 4.1]{AGGP} by adapting some crucial estimates to the three-dimensional case. 
\subsubsection{Uniqueness and continuous dependence estimate}
Let us consider two weak solutions $\ph _1$ and $\ph _2$
satisfying \eqref{wreg:1}--\eqref{wreg:2} 
and originating from two
initial data $\ph _{0}^1$ and $\ph _{0}^{2}$, where possibly $
(\varphi_{0}^1)_\Omega\neq (\varphi_0^{2})_\Omega$. Setting 
\begin{align*}
    \ph =\ph _1-\ph_2,
    \quad 
    \mu =-\ck\ast \ph + F^{\prime }(\ph _1)-F^{\prime }(\ph _2),
\end{align*}
we have
\begin{equation}
\langle \dt\ph ,v\rangle 
- \iO \ph \v\cdot \nabla v
+ \iO \nabla \mu \cdot \nabla v
=0,\quad \forall \, v \in V,\ \text{a.e. in }%
(0,T).  \label{weaky}
\end{equation}%
Taking $v=\mathcal{N}(\ph -\varphi_\Omega)$, we find
\begin{equation*}
    \frac{1}{2}\frac d {dt}
    \left\Vert \ph -\varphi_{\Omega}\right\Vert
    _{\ast }^{2}+
    \iO \ph\v \cdot \nabla  \mathcal{N}\left(\ph -\varphi_{\Omega}%
    \right)
    + \iO \mu (\ph -\varphi_{\Omega})
    =0.
\end{equation*}%
By Young's
inequality, arguing as in \cite[(4.13)]{AGGP}, we have
\begin{align}
    &  \non \iO \mu (\ph -\varphi_{\Omega}) \geq
     - \iO (\ck\ast \ph )(\ph -\varphi_{\Omega})
    + \alpha \Vert \ph \Vert ^{2}-
    \iO \left(F^{\prime }\left(\ph
    _1\right)-F^{\prime }\left(\ph _2\right)\right) \varphi_{\Omega}  
    \\ \non 
    & \quad =
     -
    \iO \nabla( \ck\ast \ph) \cdot  \nabla \mathcal{N}\left(\ph -\varphi_{\Omega}\right)
    + \alpha
    \Vert \ph \Vert ^{2}-
    \iO \left(F^{\prime }\left(\ph
    _1\right)-F^{\prime }\left(\ph _2\right)\right) \varphi_{\Omega} 
    \\ \non 
    & \quad \geq 
    -
    \Vert \ck\Vert _{W^{1,1}(B_{M})}\Vert \ph \Vert\left\Vert \ph -\varphi_{\Omega}\right\Vert _{\ast }
    + \alpha \Vert
    \ph \Vert ^{2}
    -\left\vert (\varphi_1)_\Omega-
    (\ph_{2})_\Omega\right\vert \left( \Vert F^{\prime }(\ph _1)\Vert
    _{1}+\Vert F^{\prime }(\ph _2)\Vert _{1}\right)  \\
    & \quad \geq \frac{3\alpha }{4}\Vert \ph \Vert ^{2}-C\Vert \ph -%
    \varphi_{\Omega}\Vert _{\ast }^{2}
    -\left\vert (\varphi_1)_\Omega-
    (\ph_{2})_\Omega\right\vert \left( \Vert F^{\prime }(\ph _1)\Vert
    _{1}+\Vert F^{\prime }(\ph _2)\Vert _{1}\right),
    \label{mu-diff}
\end{align}%
where $B_M$ is a sufficiently large ball containing $\overline{\Omega}$.
Concerning the convective term, by Sobolev--Gagliardo--Nirenberg's inequality, we obtain
\begin{align}
    \non \left| \iO \ph  \v \cdot \nabla \mathcal{N}\left(\ph -\varphi_{\Omega}
    \right)  \right|& \leq
    \Vert \v\Vert _6\Vert \ph \Vert
    \left\Vert \nabla \mathcal{N}\left(\ph -\varphi_{\Omega}\right)\right\Vert_3
    \\&\leq \non 
    C\Vert \v\Vert _6\Vert \ph \Vert 
    \left\Vert \nabla \mathcal{N}\left(\ph -\varphi_{\Omega}\right)\right\Vert ^{%
    \frac{1}{2}}
    \left\Vert \ph -\varphi_{\Omega}\right\Vert ^{\frac{1}{2}}\\
    & \non  \leq \frac{\alpha }{8}\Vert \ph \Vert ^{2}+C\Vert
    \v\Vert _6^{2}
    \left\Vert \ph -\varphi_{\Omega}\right\Vert _{\ast } \left(\Vert \ph \Vert +C\left|\varphi_{\Omega}\right| \right)
    \\
    & \leq \frac{\alpha }{4}\Vert \ph \Vert ^{2}
    +C\Vert
    \v\Vert _6^{4}
    \left\Vert \ph -\varphi_{\Omega} \right\Vert
    _{\ast }^{2}+C\left\vert \varphi_{\Omega}\right\vert ^{2}.
    \label{U:weak-u}
\end{align}%
Then, recalling the conservation of mass, i.e., $\ph ^{i}_\Omega(t)=%
(\varphi_{0}^{i})_\Omega$ for all $t\in \lbrack 0,T]$ and $i=1,2$, we are
led to
\begin{equation*}
    \frac d {dt}\Vert \ph \Vert _{\Vp}^{2}+\alpha \Vert \ph \Vert ^{2}\leq C\left( 1+\Vert
    \v\Vert _6^{4}\right) \Vert \ph \Vert _{\Vp}^{2}+\Lambda \left\vert \varphi_{\Omega}(0)\right\vert
    +C\left\vert \varphi_{\Omega}(0)\right\vert ^{2},
\end{equation*}%
where $\Lambda =2 (\left\Vert F^{\prime }\left(\ph _1 \right)\right\Vert _{1}+\left\Vert F^{\prime }\left(\ph_{2}\right) \right\Vert _{1})$. 
Therefore, an application of Gronwall's Lemma implies \eqref{cd}, which, in particular, entails the
uniqueness of the weak solutions $\ph$. \last{Concerning the uniqueness for the corresponding chemical potential $\mu$, it readily follows upon noticing that $\ph$ is uniquely determined in $L^2(Q)$, hence almost everywhere in $Q$, along with the fact that $F'$ is single valued.}
\medskip
\subsubsection{Existence of weak solutions}
The proof of the existence of weak solutions can be proven exactly as in \cite[Thm.4.1]{AGGP}, since all the estimates do not depend on the dimension of the domain and are thus valid also for the three-dimensional case.
\subsubsection{Existence of strong solutions: parts (i)-(ii)}
To derive the existence of strong solutions, we \rev{straightforwardly} adapt the proof of \cite[Thm.4.1]{AGGP}. 
Let us consider a
sequence of velocity fields $\{ \v^{k} \}\subset C_0^\infty(
(0,T);{\CCC}_{0,\sigma }^{\infty }(\Omega))$ such that $\v%
^{k}\rightarrow \v$ strongly in $L^{4}(0,T;\textbf{L}_\sigma^{6})$. For
any $k\in \mathbb{N}$, we introduce the Lipschitz continuous truncation $h_{k}$ given by
\begin{equation*}
    h_{k}:\mathbb{R}\rightarrow \mathbb{R},
    \quad 
    h_{k}(s)=%
    \begin{cases}
    -1+\frac{1}{k},\quad & s<-1+\frac{1}{k}, \\
    s,\quad & -1+\frac{1}{k} \leq s \leq 1-\frac{1}{k}, \\
    1-\frac{1}{k},\quad & s>1-\frac{1}{k},
    \end{cases}%
\end{equation*}
and set $\ph _{0}^{k}:=h_{k}(\ph _{0})$. It readily follows that
\begin{align*}
    \text{$\ph _{0}^{k}\in V\cap
L^{\infty }(\Omega )$ and  that $\nabla \ph _{0}^{k}=\nabla \ph _{0}
\chi _{\lbrack -1+\frac{1}{k},1-\frac{1}{k}]}(\ph _{0})$ almost everywhere
in $\Omega $,}
\end{align*}
where $\chi _{A}(\cdot )$ denotes the indicator function \an{of a measurable} set
$A$. By definition, we have
\begin{equation}
    \left|\ph _{0}^{k}\right| 
    \leq \left|\ph _{0}\right|,\quad \left|\nabla \ph _{0}^{k}\right|
    \leq \left|\nabla
    \ph _{0}\right|\quad \text{a.e. in }\Omega.
    \label{stamp}
\end{equation}%
\rev{Indeed, observe that, for a.a. $ (x,t) \in Q$, if $\varphi_0(x,t)<-1+\frac 1 k<0$, then $\vert \varphi_0(x,t)\vert=-\varphi_0(x,t)>1-\frac1k=\vert\varphi_0^k( x,t)\vert$. Estimates \eqref{stamp} entail} $\ph _{0}^{k}\rightarrow \ph _{0}$ strongly in $%
V$ as well as $\left\vert
(\ph _{0}^{k})_\Omega\right\vert \rightarrow \left\vert (\ph
_{0})_\Omega\right\vert $ as $k\rightarrow \infty$. Then, there exist positive constants $\varpi >0$ and $\overline{k}>0$ such
that
\begin{equation}
\left\vert (\ph _{0}^{k})_\Omega \right\vert \leq 1-\varpi ,\quad \forall
\,k>\overline{k}.  \label{kk}
\end{equation}%
We now notice that \cite[Thm. A.1]{AGGP} can be slightly modified to be valid in the three dimensional-case. In particular, it can be easily seen that there exists a sequence of functions $%
\{(\ph ^{k},\mu ^{k})\}$ satisfying
\begin{align}
    \label{Reg-Appr-nch:1}
    \ph ^{k}& \in L^{\infty }(0,T;V\cap L^{\infty }(\Omega)):
    \quad 
    \supess_{t\in \lbrack 0,T]}\Vert \ph ^{k}(t)\Vert_{\infty}\leq 1-\delta _{k},
    \\
    \label{Reg-Appr-nch:2}
    \ph ^{k}& \in L^{q}(0,T;W^{1,p}(\Omega )),\quad q=\frac{4p}{3(p-2)},
    \quad
    \forall \,p\in (2,\infty ), 
    \\
    \label{Reg-Appr-nch:3}
    \dt\ph ^{k}& \in L^{\infty }(0,T;\Vp)\cap L^{2}(0,T;H), 
    \\
 \label{Reg-Appr-nch:4}
    \mu ^{k}& \in H^{1}(0,T;H) \cap \C0 V\cap L^{2}(0,T;W),
\end{align}%
where $\delta _{k}\in (0,1)$ depends on $k$. The solutions satisfy
\begin{equation}
    \dt\ph ^{k}+\nabla \ph ^{k}  \cdot \v^{k} - \Delta \mu ^{k} =0,
    \quad 
    \mu ^{k}=-\ck\ast \ph ^{k} + F^{\prime } (\ph ^{k})
    \quad
    \text{ in $Q$.}
    \label{K-nCH}
\end{equation}%
In addition, $\partial _{\mathbf{n}}\mu ^{k}=0$ almost everywhere on $\Sigma$ and $\ph ^{k}(0)=\ph _{0}^{k}$ almost
everywhere in $\Omega $. Notice that the differences compared to \cite[Thm. A.1]{AGGP} are only related to the regularity \eqref{Reg-Appr-nch:2}, which comes directly from the fact that (see, e.g., \cite[(A.27), \signo{(4.48)}]{AGGP}),
\begin{align}
\Vert \nabla \ph^k\Vert _{p}\leq C^{\frac 1 p}\left( 1+\Vert
\nabla \mu ^k\Vert _{p}\right)\quad\forall p\geq 2,\label{mu6}
\end{align}%
together with the interpolation inequality 
\begin{align*}
    \Vert \nabla\mu^k\Vert_{L^{q}(0,T;\textbf{L}^{p})}
    \leq 
    C 
    \Vert \nabla\mu^k\Vert _{L^{\infty}(0,T;\textbf{H})}
    \Vert \mu^k\Vert _{L^{2}(0,T;{W})}, 
\end{align*}
where $q=\frac{4p}{3(p-2)}$ and $p\in (2,\infty )$  and $C>0$ is independent of $k$. 
Indeed, thanks to the recent result in \cite{P} (see \cite[Corollary 4.5, Remark 4.7]{P}), the strict separation property for $\ph^k$ entails the existence of  $\delta_k>0$ such that 
\begin{equation}
    \supess_{t\in \lbrack 0,T]}\Vert \ph ^{k}(t)\Vert _{\infty}\leq 1-\delta _{k}.
\label{sep}
\end{equation}
Now it is immediate to show (see also \cite[(4.22), (4.24), (4.27)]{AGGP}) that 
\begin{align}
    &
    \intQ |\nabla \mu ^{k}|^2 
    + \intQ |\nabla \ph ^{k}|^2
    + \ioT \Vert \dt\ph ^{k}\Vert _{\Vp}^{2}
    \leq 
    C(1+T)
    + C \intQ|\v^{k}|^{2}
    \leq C,
    \label{nablamu-L2}
\end{align}
where $C>0$ does not depend on $k$, since $\v^k\to\v$ strongly in $L^4(0,T;\textbf{L}_\sigma^6)$.
Observe that the regularity of the approximated solutions $\{(\ph
^{k},\mu ^{k})\}$ in \eqref{Reg-Appr-nch:1} allows us to compute the time and the
spatial derivatives of the second equation in \eqref{K-nCH}, which gives
\begin{equation}
    \dt\mu ^{k}=-\ck\ast \dt\ph ^{k} + F^{\prime \prime }(\ph ^{k})\dt\ph
    ^{k},\quad \nabla \mu ^{k}=-\nabla \ck\ast \ph ^{k} + F^{\prime \prime
    }(\ph ^{k})\nabla \ph ^{k}\quad \text{ in }%
    Q.  \label{mut}
\end{equation}%
In addition, the map $t\mapsto \Vert \nabla \mu (t)\Vert^{2}$ belongs to $AC([0,T])$ and the chain rule $\frac12\frac{d}{dt} \Vert \nabla \mu ^{k}\Vert ^{2}= \iO \dt\mu
\Delta \mu$ holds almost everywhere in $(0,T)$. Thus, testing  the first equation in 
\eqref{K-nCH} by $\dt\mu ^{k}$, integrating over $\Omega $,
and exploiting \eqref{mut}, we obtain
\begin{equation}
    \frac{1}{2}\frac d {dt}\Vert \nabla \mu ^{k}\Vert ^{2}
    +\int_{\Omega }F^{\prime \prime }(\ph ^{k})|\partial _{t}\ph ^{k}|^{2}\,
    =
    - \int_{\Omega } \nabla \ph ^{k} \cdot \v^{k}   \,\dt\mu ^{k}
    + \int_{\Omega }\ck\ast \dt\ph  ^{k}\,\dt\ph ^{k}.
    \label{mu7}
\end{equation}%
We rewrite the key term $\iO \v^{k}\cdot \nabla \ph ^{k} \partial_{t}\mu ^{k}$ as in the proof of \cite[Thm. 4.1]{AGGP}. By using \eqref{mut} and the properties of $\v^{k}$
in $C_0^\infty((0,T);{\CCC}_{0,\sigma }^{\infty }(\Omega))$, we
observe that
\begin{align}
    \non \int_{\Omega }\nabla \ph ^{k} \cdot \v^{k}  \,\dt\mu ^{k}\,%
    & =
    -\int_{\Omega }\left(\nabla \ph^{k} \cdot \v^{k}\right) \ck\ast
    \dt\ph ^{k}
    +\int_{\Omega }\left(\nabla \ph^{k} \cdot \v^{k}\right) F^{\prime \prime }(\ph ^{k})\,\dt\ph ^{k}
    \\ \non 
    & =
    -\int_{\Omega }\left(  \nabla \ph^{k} \cdot \v^{k}\right) \ck\ast \dt\ph ^{k}
    +\int_{\Omega } \nabla \left( F^{\prime }(\ph^{k})\right)\cdot\v^{k} \dt\ph ^{k}
     \\ \non 
    & =
     -\int_{\Omega }\left(  \nabla \ph ^{k} \cdot \v^{k}\right)
    \ck\ast \dt\ph ^{k}
    +\int_{\Omega }\left( \nabla \mu^{k} \cdot \v^{k}\right) \partial
    _{t}\ph ^{k}
    \\ \non    
    & \quad 
     \rev{+}\int_{\Omega }\left(  \big(%
    \nabla \ck\ast \ph ^{k}\big) \cdot \v^{k} \right) \dt\ph ^{k}
    \\ \non 
    & =
    \int_{\Omega }\left( \nabla (\ck\ast \partial
    _{t}\ph ^{k}) \cdot \v^k\right) \ph ^{k}\,
    +\int_{\Omega }\left( \nabla \mu ^{k} \cdot \v^k\right) \partial
    _{t}\ph ^{k}\,
    \\ 
    & \quad 
    \rev{+}\int_{\Omega }\left(  \big(%
    \nabla \ck\ast \ph ^{k}\big) \cdot \v^k\right) \dt\ph ^{k}\, .
    \label{drift}
\end{align}%
By exploiting the uniform $L^{\infty}$-bound of $\ph ^{k}$, we have, by standard Young's inequality for convolutions,
\signo{\begin{align}
    & \non 
    \left\vert \int_{\Omega }\left(  \big(\nabla \ck\ast \ph^{k}\big) \cdot\v^{k} \right) \dt\ph ^{k}\,\right\vert 
    \leq 
    \Vert \nabla \ck\ast \ph ^{k}\Vert_{\infty} \Vert \v^{k}\Vert \Vert \dt\ph^{k}\Vert 
    \\  
    & \quad  \leq  
    \Vert \ck\Vert _{W^{1,1}(%
    B_M)}\Vert \ph ^{k}\Vert _{\infty} \Vert \v^{k}\Vert\Vert \partial
    _{t}\ph ^{k}\Vert
    \leq \frac{\alpha}{8}\Vert \dt\ph ^{k}\Vert ^{2}
    + C\Vert \v^{k}\Vert ^{2}.
    \label{I2}
\end{align}}%
Similarly, we also find
\signo{\begin{align}
    & \non \left\vert \int_{\Omega }\left(\nabla (\ck\ast \partial
    _{t}\ph ^{k}) \cdot \v^{k} \right) \ph ^{k}\,\right\vert 
    \leq 
     \Vert \nabla \ck\ast \dt\ph ^{k}\Vert
    \Vert \v^{k}\Vert\Vert \ph ^{k}\Vert _{\infty} 
    \\ & \quad 
    \leq\Vert \ck\Vert _{W^{1,1}(%
    B_M)}\Vert \dt\ph ^{k}\Vert\Vert \v^{k}\Vert \Vert
    \ph ^{k}\Vert _{\infty} 
    \leq \frac{\alpha }{8}\Vert \dt\ph ^{k}\Vert ^{2}+C\Vert \v^{k}\Vert ^{2}.
    \label{I3}
\end{align}}%
To bound the third term on the right-hand side in \eqref{drift}, we
need a preliminary estimate of the $V$-norm of $\nabla \mu ^{k}$. To
this end, let us first observe from \eqref{K-nCH} that $\mu^k - (\mu^k)_\Omega=
\mathcal{N} (\dt\ph^k + \nabla \ph^k\cdot\v^k )$ noticing that $(\dt\ph^k + \nabla \ph^k\cdot\v^k )_\Omega=0$. \rev{The inclusion of $(\mu^k)_\Omega$ is necessary, as the mapping $\cal N$ is defined on elements of $V^*$ possessing a zero mean value. Thus through the application of $\mathcal{N}$ on $\dt\ph^k + \nabla \ph^k\cdot\v^k\in V_{(0)}^*$, we can only retrieve the $L^2(\Omega)$-projection of $\mu^k$ over the space with zero integral mean.} Then, we find, by elliptic regularity,
\begin{equation}
\Vert \nabla\mu ^{k}\Vert _{V}\leq C\left( \Vert \partial
_{t}\ph ^{k}\Vert +\Vert \nabla \ph^k\cdot\v^k\Vert \right) .  \label{muH1}
\end{equation}%
In order to estimate the second term on the right-hand side in \eqref{muH1},
we deduce from the second in \eqref{mut} that
\begin{equation}
    \nabla \ph^k\cdot\v^k=\frac{1}{F^{\prime \prime }(\ph ^{k})}%
    \left(  \nabla \mu ^{k} \cdot\v^k+(\nabla \ck\ast
    \ph ^{k}) \cdot\v^k \right) \quad \text{in }Q.  
    \label{convec}
\end{equation}%
By the strict convexity of $F$, we notice that $F^{\prime \prime
}(s)^{-1}\leq \alpha ^{-1}$ for any $s\in (-1,1)$. Thus, using Sobolev--Gagliardo--Nirenbeg's inequality and the uniform $L^{\infty }$\signo{-}bound of $\ph ^{k}$, we obtain
\begin{align}
\non \Vert\nabla \ph^k\cdot\v^k\Vert & \leq
\frac{2}{\alpha }\left( \Vert \nabla \mu ^{k}\cdot\v^k\Vert
+\Vert \big(\nabla \ck\ast \ph ^{k}\big) \cdot\v^k %
\Vert \right) \\ \non 
& \leq C\Vert \nabla \mu
^{k}\Vert _{{3}}\Vert \v^{k}\Vert _{{6}}
    +C \Vert \nabla \ck\ast \ph ^{k}\Vert _{\infty}\Vert \v^{k}\Vert 
    \\ \non 
    & 
    \leq 
    C \Vert \nabla \mu^{k}\Vert ^{\frac{1}{2}}\Vert \nabla \mu ^{k}\Vert_{V}^{\frac{1}{2}} \Vert \v^{k}\Vert _{6}
    + C \Vert \ck\Vert _{W^{1,1}(B_M)}\Vert \ph ^{k}\Vert _{\infty}\Vert \v^{k}\Vert
    \\ & 
    \leq 
    C \Vert \nabla \mu^{k}\Vert ^{\frac{1}{2}}\Vert \nabla \mu ^{k}\Vert_{V}^{\frac{1}{2}} \Vert \v^{k}\Vert _{{6}}
    + C \Vert \v^{k}\Vert .
\label{convecL2}
\end{align}%
Then, by exploiting \eqref{muH1} and \eqref{convecL2}, we infer that 
\begin{equation}
    \Vert \nabla \mu ^{k}\Vert _{V}\leq C\left( \Vert \partial
    _{t}\ph ^{k}\Vert +\Vert \nabla \mu ^{k}\Vert\Vert \v^{k}\Vert
    _{6}^{2} +\Vert
    \v^{k}\Vert \right) . 
    \label{muH1-f}
\end{equation}%
Now, again by Sobolev--Gagliardo--Nirenberg's inequality and \eqref{muH1-f}, we find
\begin{align}
    \non \left\vert \int_{\Omega }\left(  \nabla \mu ^{k}\cdot\v^{k}\right)
    \dt\ph ^{k}\,\right\vert 
    & \leq \Vert \nabla \mu ^{k}\Vert
    ^{\frac{1}{2}}\Vert \nabla \mu ^{k}\Vert _{V}^{%
    \frac{1}{2}}\Vert \v^{k}\Vert _{6}\Vert \dt\ph ^{k}\Vert
    \\ \non 
    & \leq 
   \Vert \nabla \mu ^{k}\Vert^{\frac{1}{2}}\Vert \dt\ph ^{k}\Vert^{\frac{3}{2}} \Vert \v^{k}\Vert _{6}
    + \Vert \nabla \mu ^{k}\Vert \Vert \dt\ph^{k} \Vert \Vert\v^{k}\Vert _{6}^{2}
    \\ \non 
    & \quad 
    +\Vert \nabla \mu ^{k}\Vert^{\frac{1}{2}} \Vert \v^{k}\Vert ^{\frac{1}{2}}\Vert\v^{k}\Vert _{6}\Vert \dt\ph ^{k}\Vert
    \\
    & \leq
    \frac{\alpha }{8}\Vert \dt\ph ^{k}\Vert ^{2}
    + C \Vert \nabla \mu^{k}\Vert ^{2}\Vert \v^{k}\Vert _{6}^{4} 
    + C \Vert \v^{k}\Vert ^{2}.
\label{convec-L2-2}
\end{align}%
Concerning the last term in \eqref{mu7}, we get
\begin{align}
    \int_{\Omega }\ck\ast \dt\ph ^{k}\,\dt\ph ^{k}
    & \leq \Vert \nabla \ck\ast \dt\ph ^{k}\Vert \Vert \nabla \mathcal{N}\dt\ph^k \Vert  
    \leq 
    \frac{\alpha }{8}\Vert \dt\ph\signo{^k} \Vert ^{2}+C\left( \Vert \v^{k}\Vert ^{2}+\Vert \nabla
    \mu ^{k}\Vert ^{2}\right) .
\label{Jterm}
\end{align}%
Indeed, from the first in \eqref{K-nCH}, it holds that
\begin{align}
    \label{dual}
    \Vert \nabla \mathcal{N}\dt\ph^k \Vert=\Vert \dt\ph^k \Vert_*=\Vert \partial_t\ph^k\Vert_{\Vp}\leq C
    \left( \Vert \v^{k}\Vert +\Vert \nabla
    \mu ^{k}\Vert\right).
\end{align}
Inserting the estimates \eqref{I2}, \eqref{I3}, \eqref{convec-L2-2} and %
\eqref{Jterm} in \eqref{mu7}, we end up with
\begin{equation}
    \frac{1}{2}\frac d {dt}\Vert \nabla \mu ^{k}\Vert
    ^{2}+\frac{\alpha }{2}\Vert \dt\ph ^{k}\Vert
    ^{2}\leq C\left( 1+\Vert \v^{k}\Vert _{6}^{4}\right) \Vert \nabla \mu ^{k}\Vert ^{2}+C\Vert
    \v^{k}\Vert ^{2}. 
    \label{mut2}
\end{equation}%
Arguing as in \last{\cite[(4.49)]{AGGP}}, we then obtain
\begin{equation}
    \Vert \nabla \mu ^{k}(0)\Vert \rightarrow \Vert -\nabla \ck\ast \ph _{0} + F^{\prime
    \prime }(\ph _{0})\nabla \ph _{0}\Vert
    \quad \text{as }k\rightarrow \infty .  
    \label{nablamu0}
\end{equation}%
Therefore, since we know that $\v^k\to \v$ strongly in $L^4(0,T;\textbf{L}^6_\sigma)$ and thus it is bounded in the same space, Gronwall's lemma then entails 
\begin{equation}
\Vert \nabla\mu^k\Vert_{L^\infty(0,T;\HHH)}+\Vert \partial_t\ph^k\Vert_{L^2(0,T;H)}\leq C,
\end{equation}%
from which, recalling that it is standard to obtain $\Vert \mu^k\Vert_V\leq C(1+\Vert \nabla\mu^k\Vert)$, from \eqref{muH1-f}, we also deduce
\begin{equation}
\Vert \mu^k\Vert_{L^\infty(0,T;V) \cap L^2(0,T;W)}\leq C,
\label{a9}
\end{equation}
uniformly in $k$. Concerning the concentration $\ph ^{k}$, we deduce from \eqref{mu6}, %
\eqref{a9} and interpolation, that
\begin{equation}
\Vert \ph ^{k}\Vert _{L^{\infty }(Q)}\leq 1,\quad \Vert
\ph ^{k}\Vert _{L^{\infty }(0,T;V)\cap L^{q}(0,T;W^{1,p}(\Omega ))}\leq C,\quad q=\frac{4p}{3(p-2)}. \label{a3}
\end{equation}%
In a similar fashion, by comparison in \eqref{K-nCH}, we are led to
\begin{equation}
\Vert F^{\prime }(\ph ^{k})\Vert _{L^{\infty }(0,T;V)\cap L^{q}(0,T;W^{1,p}(\Omega ))}\leq C,
\label{a4}
\end{equation}%
with the same $q$ as above.
Furthermore, recalling \eqref{dual}, we obtain from $\v \in
L^{4}(0,T;\textbf{L}_{\sigma }^{6}(\Omega ))$ and \eqref{a9} that
\begin{align}
    \Vert \dt\ph ^{k}\Vert _{L^{4}(0,T;\Vp)}\leq C.  
    \label{a5}
\end{align}
Exploiting the above uniform estimates, by a
standard compactness argument, we pass to the limit in a suitable weak form of \eqref{K-nCH} as $k\rightarrow \infty $
obtaining that the limit function $\ph $ is a strong solution 
to \eqref{eq:1}-\eqref{eq:4}. In particular, \eqref{eq:1} and \eqref{eq:2} hold almost
everywhere in $Q$, \signo{and} \eqref{eq:3} holds almost
everywhere on $\Sigma$. 
Finally, if we also assume $\v\in L^{\infty }(0,T;\textbf{H}_\sigma)$ \signo{(cf. \eqref{dual})} it
is easily seen that $\dt\ph \in L^{\infty }(0,T;\Vp)$. This concludes the proof of part $(i)$. \\

\noindent 
Concerning part $(ii)$, as already noticed for the approximating sequence, it is a consequence of \cite[Corollary 4.5, Remark 4.7]{P}. Indeed, if the initial datum is strictly separated from pure phases, then the strict separation property for $\ph$ holds as well, thanks to the regularity of the solution $\ph$ and the divergence-free property of the advective vector field $\v$, so that there exists $\delta>0$, depending on $\varphi_0$, such that 
\begin{equation}
    \supess_{t\in \lbrack 0,T]}\Vert \ph (t)\Vert _{\infty}\leq 1-\delta .
    \label{sepa}
\end{equation}
Arguing by means of the difference  quotients (see \cite[Thm. 4.1, part (ii)]{AGGP}), it is then easy to show  that $\dt\mu \in
L^{2}(0,T;H)$. In turn, thanks to $\mu \in
L^{2}(0,T;W)$, we also obtain $\mu \in \C0 V$ \signo{via compact embeddings}.%
\medskip
\subsubsection{Extra regularity results for separated initial data: part (iii)}
We perform the proof in the same spirit of \cite[Lemma 2]{Frigeri}. Let $\kk\in [-1+\delta,1-\delta]$ and $\eta = \eta (x, t) \in [0, 1]$ be a continuous piecewise-smooth function which is supported on the space-time cylinders $Q_{t_0,t_0+\tau} (\rho) := B_{\rho} (x_0) \times
(t_0, t_0 + \tau )$, where $B_\rho (x_0)$ denotes the ball centered at $x_0$ of radius $\rho> 0$, and $\tau>0$ is given. According to the type of H\"{o}lder regularity we consider (interior or boundary regularity), we set $x_0\in \Omega$ or $x_0\in\Gamma$ and then by standard compactness arguments we gain the validity for any ball of $\overline{\Omega}$. We thus test \eqref{weak-nCH-2:1} by $\eta^2\ph^+_\kk$, where $\ph^+_\kk:= \max \{0, \ph - \kk\}$, integrate the resulting identity over $Q_{t_0,t} :=   \Omega\times(t_0, t)$, where $0 \leq t_0 < t < t_0 + \tau \leq T$ , to infer that, \last{for any $s\leq t$},
{\begin{align}
   \nonumber&\frac 12 \Vert \eta\varphi_\kk^+ (s)\Vert^2
    + \int_{Q_{t_0,s}}F''(\ph)\nabla\varphi_\kk^+\cdot \nabla(\eta^2\varphi_\kk^+)
   =\frac 12 \Vert {\eta}\varphi_\kk^+(t_0)\Vert^2+\int_{Q_{t_0,s}}(\nabla \ck\ast \ph)\cdot\nabla(\eta^2\varphi_\kk^+)
   \nonumber\\&\quad
   +\int_{Q_{t_0,s}}\varphi_\kk^+\v\cdot \nabla(\varphi_\kk^+\eta^2)+{\int_{Q_{t_0,s}}(\ph_\kappa^+)^2\eta\partial_t\eta}.
    \label{esta}
\end{align}}%
{\last{In the above identity,} we exploited some properties of the {positive part function. For instance, we used that} 
\begin{align*}
    \int_{Q_{t_0,s}}F''(\ph)\nabla\varphi\cdot \nabla(\eta^2\varphi_\kk^+)
    &=\int_{\{(x,r)\in Q_{t_0,s}:\ \ph(x,r)\geq \kk\}}F''(\ph)\nabla\varphi\cdot \nabla(\eta^2{(\varphi-\kk)})
    \\&
    =\int_{Q_{t_0,s}}F''(\ph)\nabla\ph_\kk^+\cdot \nabla(\eta^2\varphi_\kk^+),
\end{align*}
and {that}
\begin{align*}    \int_{Q_{t_0,s}}\partial_t\ph \,\ph_\kappa^+\eta^2
    &=\int_{\{(x,r)\in Q_{t_0,s}:\ \ph(x,r)\geq \kk\}}{\partial_t \ph(\varphi-\kk)}\eta^2
    =\int_{Q_{t_0,s}}\frac 1 2 \partial_t(\ph_\kk^+)^2\eta^2\\&
    ={\frac 12}\norma{{\eta}\ph_\kk^+ (s) }^2-{\frac 12}\norma{{\eta}\ph_\kk^+(t_0)}^2-\int_{Q_{t_0,s}}(\ph_\kappa^+)^2\eta\partial_t\eta.
\end{align*}
}%
Now, we point out  the identity
$$
    \nabla\varphi_\kk^+\cdot \nabla(\eta^2\varphi_\kk^+)
    =\vert \nabla(\varphi_\kk^+\eta)\vert^2-\vert \nabla\eta\varphi_\kk^+\vert^2,
$$
{{whence} \eqref{esta} entails, {using again \ref{h2}, that} 
\begin{align}
   \nonumber&\frac 12 \sup_{s\in[t_0,t]}\Vert {\eta}\varphi_\kk^+(s)\Vert^2
    + \alpha\int_{Q_{t_0,t}}\vert\nabla(\varphi_\kk^+ \eta)\vert^2
    \\&\quad
   \leq \frac 12 \Vert {\eta}\varphi_\kk^+(t_0)\Vert^2+\int_{Q_{t_0,t}}\left\vert(\nabla \ck\ast \ph)\cdot\nabla(\eta^2\varphi_\kk^+)\right\vert\nonumber\\&\qquad+\sup_{s\in[t_0,t]}\left\vert\int_{Q_{t_0,s}}\varphi_\kk^+\v\cdot \nabla(\varphi_\kk^+\eta^2)\right\vert+{\int_{Q_{t_0,t}}(\ph_\kappa^+)^2\vert\eta\partial_t\eta\vert+\int_{Q_{t_0,t}}F''(\ph)\vert\nabla\eta\varphi_\kk^+\vert^2}.
    \label{est2}
\end{align}
}%
{Using the separation property \eqref{sepa}, it holds that} $\Vert F''(\ph)\Vert_{L^\infty(Q)}\leq C$. Thus, 
$$
    \int_{Q_{t_0,t}}F''(\ph)\vert\nabla\eta\varphi_\kk^+\vert^2\leq C\int_{Q_{t_0,t}}\vert\nabla\eta\vert^2(\varphi_\kk^+)^2.
 $$
Moreover, we {readily} deduce that
\begin{align*}
    & \int_{Q_{t_0,t}}\left\vert(\nabla \ck\ast \ph)\cdot\nabla(\eta\varphi_\kk^+{\eta})\right\vert
    \\&\quad 
    \leq 
    \Vert \nabla \ck\ast \ph\Vert_{L^\infty(Q)}\left(\int_{Q_{t_0,t}}\vert \nabla\eta\vert\varphi_\kk^+\eta +\int_{Q_{t_0,t}}\eta\vert\nabla(\varphi_\kk^+\eta)\vert \right)\\&\quad 
    \leq \frac{\alpha}{4}\int_{Q_{t_0,t}}\vert\nabla(\eta\varphi_\kk^+)\vert^2+C\int_{Q_{t_0,t}}\eta^2+C\int_{Q_{t_0,t}}\vert \nabla\eta\vert^2(\varphi_\kk^+)^2.
\end{align*}
In conclusion, we observe that
$$
  {\int_{Q_{t_0,s}}\varphi_\kk^+\v\cdot \nabla(\varphi_\kk^+\eta^2)=-\int_{Q_{t_0,s}}
    \nabla \Big(\frac{(\varphi_\kk^+)^2}{2}\Big)
    \cdot \v \eta^2=\int_{Q_{t_0,s}}\eta(\varphi_\kk^+)^2\v\cdot \nabla\eta} ,
$$
and, by \last{the} assumption $\Vert\v\Vert_{L^4(0,T;\textbf{L}^6)}\leq C$ and standard inequalities, we get
{\begin{align*}
   &\sup_{s\in[t_0,t]}\left\vert\int_{Q_{t_0,s}}\varphi_\kk^+\v\cdot \nabla(\varphi_\kk^+\eta^2)\right\vert\\&\quad\leq 
    \int_{Q_{t_0,t}}\left\vert\eta(\varphi_\kk^+)^2\v\cdot \nabla\eta\right\vert
   \\&\quad \leq \Vert \eta\varphi_\kk^+\Vert_{    {L^4(t_0,t;\Lx3)}}\Vert \varphi_\kk^+\nabla\eta\Vert_{L^2({Q_{t_0,t}})}\Vert \v\Vert_{L^4(0,T;\textbf{L}^6)}
    \\&\quad 
    \leq {C\left[\int_{{t_0}}^{t}\Vert \eta\varphi_\kk^+\Vert^2\left(\Vert \eta\varphi_\kk^+\Vert^2
    +\Vert \nabla(\eta\varphi_\kk^+)\Vert^2\right)\right]^{\frac 1 4}\Vert \varphi_\kk^+\nabla\eta\Vert_{L^2(Q_{t_0,t})}}
    \\&\quad \leq
    \frac{\alpha}{4}\int_{Q_{t_0,t}}\vert\nabla(\eta\varphi_\kk^+)\vert^2+\frac 1 4\sup_{s\in[t_0,t]}\Vert \eta\varphi_\kk^+(s)\Vert^2+C\int_{Q_{t_0,t}}\vert \nabla\eta\vert^2(\varphi_\kk^+)^2.
\end{align*}}
Summarizing, putting everything together in \eqref{est2}, we deduce
\begin{align}
    \nonumber&\frac 1 4\sup_{s\in[t_0,t]}\Vert \eta\varphi_\kk^+(s)\Vert^2
    +\frac \alpha 2\int_{Q_{t_0,t}}\vert\nabla(\eta\varphi_\kk^+)\vert^2
    \\&\quad 
    \leq
    \frac 1 2\Vert \eta\varphi_\kk^+(t_0)\Vert^2
    +C\int_{Q_{t_0,t}}\vert \nabla\eta\vert^2(\varphi_\kk^+)^2+C\int_{Q_{t_0,t}}\eta^2+{\int_{Q_{t_0,t}}(\ph_\kappa^+)^2\vert\eta\partial_t\eta\vert}.
    \label{esta2}
\end{align}
Arguing in a similar way, inequality \eqref{esta2} also holds with $\ph$ when replaced by $-\ph$ leading us to consider {$\varphi_\kk^-=(-\ph-\kk)^+$}. In particular, for any fixed $\varepsilon>0$ such inequalities imply that $\ph$ is an element of
$\mathfrak{B}_2(\ov Q,1-\delta, \gamma_, \omega, \varepsilon,\chi)$ in the sense of \cite[Ch. II, Sec. 7]{Lady}, for some $\gamma,\omega,\chi> 0$ possibly depending on $\varepsilon$ (cf., in particular, the inequalities in \cite[\lasts{Ch.} V, (1.12)-(1.13)]{Lady}). Therefore, on account of \cite[Ch. V, Thm 1.1]{Lady}, the \Holder\ regularity \eqref{holderphi} holds. Then, by the regularity of $F$ and by \eqref{sepa}, we immediately deduce the same result \eqref{holder} for the chemical potential $\mu$, concluding the proof.

\subsubsection{$L^4_tW^{1,6}_x$ regularity on {$\ph$} for separated initial data: part (iv)}
\last{Let us now prove point $(iv)$ of Theorem  \ref{ExistCahn}.} The proof can be carried out by carefully combining a version of a maximal regularity theorem in \cite{Pruss} with a \signo{bootstrap} argument, based upon the H\"{o}lder regularity in \eqref{holderphi}-\eqref{holder}.
Being the equations \eqref{eq:1}-\eqref{eq:2} satisfied by our unique solution $(\ph,\mu)$ almost everywhere in $Q$, and being $\partial_t\mu\in L^2(Q)$ as a consequence of point $(ii)$, we may rewrite \Sys\ in the form
\begin{alignat*}{2}
    & \partial_t \mu -F''(\ph)\Delta\mu=-\ck\ast \partial_t\ph -\Fsec(\ph)\nabla\ph\cdot\v=:f && \quad\text{in $Q$},\\
    & \partial_\textbf{n}\mu=0 && \quad\text{on $\Sigma$},\\
    & \mu(0)=\mu_0 && \quad \text{in $\Omega$},
\end{alignat*}
where $\mu_0$ is defined through \eqref{def:muz}.
From the regularity \eqref{holderphi}-\eqref{holder}, we obtain, for some $\beta \in (0,1)$ as in the statement, that
\begin{align}
\ph,\mu \in C^{\beta,\frac \beta 2}(\ov Q)\hookrightarrow L^\infty(0,T;C^{\beta }(\overline{\Omega})).
\label{regulmu}
\end{align}
Moreover, owing to \eqref{reg:4}, and using the limit version of  \eqref{mu6} (which is also valid in the limit as $k\to\infty$), we get
\begin{align}
    \Vert\ph\Vert_{W^{1,p}}\leq C^{{\frac 1 p}}(1+\Vert\mu\Vert_{W^{1,p}}),
    \quad p \geq 2.
    \label{pp}
\end{align}
We recall the following embedding (see, e.g., \cite[\signo{p. 318}]{Frig}): \last{for any $g\in C^\beta(\overline{\Omega})$,} 
\begin{align}
    \Vert g\Vert_{W^{\vartheta,s}}\leq C\Vert g\Vert_{C^\beta(\overline{\Omega})}, \quad  \vartheta\in(0,\beta),\quad  s\in(1,\infty).
    \label{emb}
\end{align}
Moreover, for any $g\in W^{\th,s}(\Omega)\cap W^{2,\gamma }(\Omega)$, with  \begin{align}
  \frac{1}{2s}=\frac{1}{q}-\frac {1} {2\gamma}
  +\frac{\th} {6},\quad \th \in (0,\beta),\quad s\in(1,\infty),\quad \gamma\geq 2,
\label{s1}
\end{align}
it holds the inequality \an{(see, e.g., \cite{BM})}
\begin{align}
    \Vert  g \Vert^4_{W^{1,q}}\leq C\Vert g\Vert_{W^{\vartheta,s}}^{2}\Vert g \Vert_{W^{2,\gamma}}^{2}.
    \label{boot}
\end{align}
Let us fix $\th\in(0,\beta)$ and assume that $\frac1 q= \frac{1}{2\gamma}-\frac \th{12}$, i.e., $q=\frac{12\gamma}{6-\gamma\th}>2\gamma\geq 4$. Then, to accomplish \eqref{s1} for any $2\leq \gamma< \frac{6}{\vartheta}$, notice that $\frac{6}{\vartheta}>6$, it suffices fixing $s=\frac 6 \th\in({6},\infty)$.
We then introduce the real sequence $\{q_n\}$  recursively by setting
\begin{align}
  4<q_1
  = \frac{12}{3-\th}\leq\ \ldots
  \leq q_n
  = \frac{12}{6-\widetilde{q}_{n}\th}\ \widetilde{q}_{n},
  \quad 
  \widetilde{q}_{n}
    :=\frac{6q_{n-1}}{6+q_{n-1}},
  \label{sequence}
\end{align} 
where the usefulness of the auxiliary sequence $\{\widetilde q_n\}$ will be clarified below.
Letting $n \to \infty$, we find that 
\begin{align*}
    q_n 
    \to \frac{6}{1-\th}>6
    \quad 
    \text{and}
    \quad
    \widetilde{q}_{n}
   \to \frac{6}{2-\th}>3.
\end{align*}
Besides, being $\th\in(0,1)$, it holds that ${4}<q_n\leq \frac{6}{1-\th}$ for every $n\geq 1$ and it is monotone increasing: whence it possesses a limit that we claim to be finite.
If that is not the case, meaning that  $q_n\to +\infty$ as $n \to \infty$, we would get a contradiction since by construction $\widetilde{q}_n\to 6$ and $q_n\to \frac{12}{1-\th}<\infty$. Therefore, the finite value $\ov q :=  \frac{6}{1-\th}$ is the limit of the sequence and we also notice that  $\widetilde{q}_{n}\to \frac{6}{2-\th}$ as $n\to\infty$. From \last{\eqref{emb}-\eqref{boot}}, along with the regularity  of $\mu$ in \last{\eqref{regulmu}}, we thus get
\begin{align}
    \Vert  \mu \Vert_{L^{4}(0,T;\Wx{1,{q_1}})}^{4}\leq C\Vert \mu\Vert_{L^\infty(0,T;\Wx{\vartheta,s})}^{2}\Vert \mu \Vert_{L^2(0,T;H^2(\Omega))}^{2}\leq C,
    \label{boot2}
\end{align}
so that, from \eqref{pp}, we deduce 
\begin{align}
    \Vert  \ph\Vert_{L^{4}(0,T;\Wx{1,{q_1}})}\leq C.
    \label{boot3b1}
\end{align}
Following the same arguments as in the proof of \cite[Lemma 5.4]{P}, being now $\ph$ bounded in $L^4(0,T;W^{1,4}(\Omega))$ since $q_1>4$, we immediately infer \eqref{p2} without extra assumptions on the initial dat\an{a $\varphi_0,\mu_0$} with respect to point $(iii)$.

\last{Furthermore}, if $\mu_0\in {\BB}_{3,2}^{1}(\Omega)$, we can apply a regularity result presented in \cite[Sec.4]{Pruss}. Namely, recalling \eqref{holderphi}, we have, with the notation of the quoted paper\last{,} 
\begin{align*}    a_{kl}&:=\delta_{kl}F''(\ph)\in C^{\beta,\frac \beta 2}(\ov Q),\quad a_k=0,\quad b_0=0,\quad b_k=\an{n}_k,\quad k,l=1,2,3\\
        f&:=-\ck\ast \partial_t\ph -\Fsec(\ph)\nabla\ph\cdot\v ,
\end{align*}
where $\an{n}_k$ is the $k$-th component of the outward unit normal \an{$\nnn$}, which is smooth by assumption, and 
$\delta_{kl}$ are the Kronecker delta symbols. Clearly, $\{a_{kl}\}_{kl}$ is uniformly elliptic since $F''\geq \alpha>0$ by \ref{h2}. Observe now that, for any $p\leq 6$, we have
\begin{align}
    \norma{\ck\ast \partial_t\ph}_{L^2(0,T;L^p(\Omega))}\leq C\norma{\ck\ast \partial_t \ph}_{L^2(0,T;V)}\leq C\norma{\partial_t\ph}_{L^2(Q)}\leq C,
\label{dtp}
\end{align}
by previous bounds.
Furthermore, we have from the separation property \eqref{SP} that $\Vert F^{''}(\ph)\Vert_{\infty}\leq C$, and so for any $n\geq2$
\begin{align}
\norma{\Fsec(\ph)\nabla\ph\cdot\v}_{L^2(0,T;L^{\widetilde{q}_n}(\Omega))}\leq C\norma{\nabla\ph}_{L^4(0,T;\textbf{L}^{q_{n-1}}(\Omega))}\norma{\v}_{L^4(0,T;\LLL_\sigma^6)},
\label{tricks}
\end{align}
where $\widetilde{q}_n$ is the $n$-th element of the sequence introduced above.
Consider $n=2$: we have  
\begin{align}
\norma{\Fsec(\ph)\nabla\ph\cdot\v}_{L^2(0,T;L^{\widetilde{q}_2}(\Omega))}\leq C\norma{\nabla\ph}_{L^4(0,T;\textbf{L}^{q_{1}}(\Omega))}\norma{\v}_{L^4(0,T;\LLL_\sigma^6)}\leq C,
\label{tricks2}
\end{align}
and thus, together with \eqref{dtp}, it yields that 
$$
\norma{f}_{L^2(0,T;L^{\widetilde{q}_2}(\Omega))}\leq C.
$$
Then, recalling that $\mu_0\in {\BB}_{3,2}^{1}(\Omega)\hookrightarrow {\BB}_{q,2}^{1}(\Omega)$, for any $q\leq 3$ 
(as long as ${\widetilde{q}_n}<3$, otherwise we would be already done), 
applying the same result as above (cf. \cite[Sec.4]{Pruss}), we get 
\begin{align}
 \Vert  \mu\Vert_{ \H1 {\Lx {\widetilde q_2}} \cap L^2(0,T;W^{2,{\widetilde{q}_2}}(\Omega))}\leq C.
	\label{boot6}\end{align}
Coming back to \eqref{boot}, holding with $\gamma={\widetilde{q}_2}\in[2,\frac{6}{\th})$, we now infer 
	\begin{align}
    \Vert  \mu \Vert_{L^{4}(0,T;\Wx{1,{q_2}})}^{4}\leq C\Vert \mu\Vert_{L^\infty(0,T;\Wx{\vartheta,s})}^{2}\Vert \mu \Vert_{L^2(0,T;W^{2,\widetilde{q}_2}(\Omega))}^{2}\leq C,
    \label{boots}
\end{align}
which then entails, by \eqref{pp}, that
\begin{align}
    \Vert  \ph\Vert_{L^{4}(0,T;\Wx{1,{q_2}})}\leq C.
    \label{boot3b}
\end{align}
Then, we are essentially ready to iterate. Namely, by the same arguments as above, we infer that
$$
\norma{f}_{L^2(0,T;L^{\widetilde{q}_3}(\Omega))}\leq C,
$$
so that by \cite[Sec.4]{Pruss}, \signo{\eqref{pp} and \eqref{boot}} (as long as ${\widetilde{q}_3}<3$, otherwise we would be already done), we get 
\begin{align}
 \Vert  \mu\Vert_{\H1 {{\Lx {\widetilde q_3}}} \cap L^2(0,T;W^{2,{\widetilde{q}_3}}(\Omega))}+\Vert  \ph\Vert_{L^{4}(0,T;\Wx{1,{q_3}})}\leq C.
	\label{boot6b3}\end{align}
Finally, by extending this \signo{bootstrap} argument, at a general step $n$ for which $q_n<6$ and $\widetilde{q}_n<3$, we obtain that
\begin{align}
 \Vert  \mu\Vert_{\H1 {\Lx {\widetilde q_n}} \cap L^2(0,T;W^{2,{\widetilde{q}_n}}(\Omega))}+\Vert  \ph\Vert_{L^{4}(0,T;\Wx{1,{q_n}})}\leq C.
	\label{boot6b}\end{align}
According to \eqref{sequence}, we have $q_n\to \frac{6}{1-\th}>6$ and $\widetilde{q}_n\to \frac{6}{2-\th}>3$, so that there exists $\overline{n}>0$ such that 
$$
    \floor*{q_{\overline{n}}}= 6,\quad \floor*{\widetilde{q}_{\ov{n}}}= 3,
$$
where $\floor{\cdot}$ denotes the floor function.
Therefore, we can iterate the arguments leading to \eqref{boot6b} starting from step $\ov{n}-1$, with the above two quantities in the resulting summability exponents, recall that $\mu_0\in\BB_{3,2}^1(\Omega)$, and thus prove \eqref{regul}.
This concludes the proof of the theorem since by comparison we immediately deduce the desired regularity on 
\an{$\partial_t\ph$} as well.
\subsubsection{{$ \signo{ L^4_tH^{2}_x\cap L^3_t W^{1,\infty}_x}$ regularity on $\ph$} for separated initial data: part (v)}
\last{We are left to prove point $(v)$ of Theorem~\ref{ExistCahn}.
Recall that we are \an{additionally} assuming $\v\in  L^4(0,T;\LLL^\infty)$ as well as $\mu_0\in \BB^1_{6,2}(\Omega)$.}
Then, by exploiting the results of the previous section, the regularity in \eqref{regul} is fulfilled {and we easily reach  }
$$
    \norma{\Fsec(\ph)\nabla\ph\cdot\v}_{L^2(0,T;L^6(\Omega))}\leq C\norma{\nabla\ph}_{L^4(0,T;\textbf{L}^6(\Omega))}\norma{\v}_{L^\an{4}(0,T;\LLL_\sigma^\an{\infty})}\leq C. 
$$
{Going back to} \eqref{dtp}, with the same notation as the previous section, we infer that
$$
    f=-\ck\ast \partial_t\ph-\Fsec(\ph)\nabla\ph\cdot\v\in L^2(0,T;L^6(\Omega)),
$$
so that, {owing to the same result as above (\cite[Sec.4]{Pruss}),} we immediately infer {that}
\begin{align}
 \Vert  \mu\Vert_{ \H1 {\Lx {6}} \cap L^2(0,T;W^{2,{6}}(\Omega))}\leq C.
	\label{boot9}\end{align}
We now recall (see{, e.g.,} \cite{BM}) that it holds the inequality 
\begin{align}
    \Vert  \mu \Vert^3_{W^{1,\infty}}\leq C\Vert \mu\Vert_{W^{\vartheta,s}}\Vert \mu \Vert_{W^{2,6}}^{2},
    \label{boota}
\end{align}
where, as in the previous section, $\theta\in(0,\beta)$ {with $\beta$ being the \Holder\ exponent in \eqref{holder},} and $s=\frac{3}{\theta}$. This {entails, using} \eqref{regulmu} and \eqref{boot9}\an{,} that
$$
    \norma{\mu}_{L^3(0,T;W^{1,\infty}(\Omega))}\leq C,
$$
so that, by letting $p\to \infty$ in \eqref{pp}, this {implies that}
$$
    \norma{\ph}_{L^3(0,T;W^{1,\infty}(\Omega))}\leq C.
$$
\an{Then, by comparison, we immediately deduce the \signo{required} regularity on 
\an{$\partial_t\ph$} as well.}

\noindent{Next, assume $\v\in L^\infty(0,T;\LLL_\sigma^4)\cap L^4(0,T;\LLL_\sigma^6)$ and $\mu_0\in \BB^{\frac{3}2}_{2,4}(\Omega)$.} In this case we have 
$$
    \norma{\Fsec(\ph)\nabla\ph\cdot\v}_{L^4(0,T;H)}\leq C\norma{\nabla\ph}_{L^4(0,T;\textbf{L}^4(\Omega))}\norma{\v}_{L^\infty(0,T;\textbf{L}^4(\Omega))}\leq C,
$$
again thanks to {\eqref{regul}.}
Furthermore, {being $\ck$} symmetric, for any $\psi\in V$,
$$
    \vert\langle \ck\ast \partial_t\ph, \psi \rangle\vert={\Big|\iO \ck\ast \partial_t\ph\psi\Big|=\Big| \iO \partial_t\ph\ck\ast \psi \Big|}\leq C\norma{\partial_t\ph}_{\an{*}}\norma{\ck\ast \psi}_V\leq C\norma{\partial_t\ph}_{\an{*}}\norma{\psi}_V,
$$
entailing \signo{that}
$$
    \norma{\ck\ast \partial_t\ph}_{L^\infty(0,T;V^*)}\leq C\norma{\partial_t\ph}_{L^\infty(0,T;V^*)}.
$$
 Therefore, \signo{by} standard interpolation, recalling the regularity \signo{of} $\ck$ and Young's inequality for convolutions,
\begin{align*}
    \norma{\ck\ast\partial_t\ph}_{L^4(0,T;H)}&\leq C\norma{\ck\ast\partial_t\ph}_{L^\infty(0,T;V^*)}^{\frac 1 2}\norma{\ck\ast\partial_t\ph}_{L^2(0,T;V)}^\frac 1 2\\&\leq C\norma{\partial_t\ph}_{L^\infty(0,T;V^*)}^{\frac 1 2}\norma{\partial_t\ph}_{L^\an{2}(0,T;H)}^{\frac 1 2}\leq C,
\end{align*}
\an{where the last constant $C>0$ appears exploiting the regularity given in part (i), since $\v \in L^\infty(0,T;\LLL^4_\sigma)\hookrightarrow L^\infty(0,T;\HHH_\sigma)$.} 
\an{We thus} have
$$
    f=-\ck\ast \partial_t\ph-\Fsec(\ph)\nabla\ph\cdot\v\in L^4(0,T;H),
$$
so that, by the result in \cite[Sec.4]{Pruss} {as above,} we immediately infer 
\begin{align}
 \Vert  \mu\Vert_{ W^{1,4}(0,T;H) \cap L^4(0,T;H^2(\Omega))}\leq C.
	\label{boot10}\end{align}
\an{Now it holds (see\signo{, e.g.,} \cite{BM}) \signo{that}
\begin{align}
    \Vert  \mu \Vert^8_{W^{1,4}}\leq C\Vert \mu\Vert_{W^{\vartheta,s}}^4\Vert \mu \Vert_{H^2}^{4},
    \label{boota2}
\end{align}
where $\signo{\th}\in(0,\beta)$ and $s=\frac{3}{\signo{\th}}$,} 
so that \an{by \eqref{regulmu} and \eqref{boot10},} $\mu\in L^8(0,T;W^{1,4}(\Omega))$ and thus, by \eqref{pp}, also $\ph\in L^8(0,T;W^{1,4}(\Omega))$. This allows to argue as in the proof \cite[Lemma 5.4]{P}, to deduce that 
$$
\norma{\ph}_{L^4(0,T;H^2(\Omega))}\leq C,
$$
concluding the proof of the theorem.

\subsection{Continuous dependence results}

In this subsection, we aim at proving the continuous dependence results stated in Theorem \ref{THM:CD} with respect to the velocity field $\v$. This will be crucial to establish some differentiability properties of the control-to-state mapping \last{associated} to the control problem that will allow us to identify the first-order optimality conditions for minimizers.

\begin{proof}[Proof of Theorem \ref{THM:CD}]

To begin with, we recall the notation of the statement and set 
\begin{align*}
	\v:= \v_1-\v_2,
	\quad 
	\ph:= \ph_1-\ph_2,
	\quad 
	\mu:= \mu_1-\mu_2,
 \quad 
 \ph_0:= \ph_0^1 - \ph_0^2,
\end{align*}
and write the system of the differences using \Sys.
This leads us to
\begin{alignat}{2}
	\label{eq:cd:1}
	& \dt \ph  {+\nabla \ph \cdot \v_1} + \nabla \ph_2 \cdot \v - \Delta \mu
	=
	0
	\quad && \text{in $Q$,}
	\\
	\label{eq:cd:2}
	& \mu =  - \ck \ast \ph + (F' (\ph_1)-F' (\ph_2)) 
	\qquad && \text{in $Q$,}
	\\
	\label{eq:cd:3}
	& \dn \mu = 
	0
	\quad && \text{on $\Sigma$,} 
	\\
	\label{eq:cd:4}
	& \ph(0)=\ph_0
	\quad && \text{in $\Omega$.}
\end{alignat}
\Accorpa\Syscd {eq:cd:1} {eq:cd:4}

\subsubsection*{\signo{First estimate: part (i)}}
We now observe that, thanks to the regularity of the two solutions under consideration, the following  estimates are actually rigorous (see also Remark \ref{REM:regula}). In particular, by \signo{the assumptions}, we have, for some $\delta>0$,
\begin{align}
\Vert \ph_i\Vert_{\infty}\leq 1-\delta,\ i=1,2,\quad \text{and} \quad \norma{\ph_2}_{L^4(0,T;{W}^{1,6}(\Omega))}\leq C,
    \label{reg}
\end{align}
thanks to {\eqref{regul}, holding for $\ph_2,$}\signo{.}


We test \eqref{eq:cd:1} by $\ph$, take the gradient of \eqref{eq:cd:2} and test it by $\nabla \ph$. 
In those computations, let us highlight that 
\begin{align*}	\nabla (F' (\ph_1)-F' (\ph_2))
	& = F'' (\ph_1) \nabla \ph_1 - F'' (\ph_2)\nabla \ph_2
	\\& 
	= (F'' (\ph_1) -F'' (\ph_2) ) \nabla \ph_2 + F'' (\ph_1)\nabla \ph.
\end{align*}
Then, adding the resulting identities and rearranging the terms, \an{we obtain}
\begin{align*}
	& \frac 12 \frac d {dt} \norma{\ph}^2
	+ \iO  F''(\ph_1)|\nabla \ph|^2
	= 
	-\iO (\nabla \ph_2 \cdot \v) \ph
	+ \iO (\nabla \ck \star  \ph )\cdot \nabla\ph
	\\ & \quad 
	- \iO  (F'' (\ph_1) -F'' (\ph_2) ) \nabla \ph_2 \cdot \nabla \ph.
\end{align*}
Here, we also account of the fact that one term vanishes as
\begin{align*}
     \iO \nabla \ph \cdot \v_1 \ph =
     \iO \nabla \Big(\frac {{\ph}^2} 2\Big) \cdot \v_1
    = 
    - \iO \Big(\frac {{\ph}^2} 2\Big)\div \v_1
    + \iG \Big(\frac {{\ph}^2} 2\Big)  \v_1 \cdot \nnn= 0.
\end{align*}
Moreover, we owe to \ref{h2} to infer that 
\begin{align*}
     \iO F''(\ph_1)|\nabla \ph|^2 \geq \alpha \norma{\nabla \ph}^2.
\end{align*}
Then, for a positive \signo{constant} $\d$ yet to be selected, we have
\begin{align*}
    -\iO (\nabla \ph_2 \cdot \v) \ph
	  & 
    =
    \iO \ph_2  \v \cdot \nabla \ph
    \leq
    \norma{\ph_2}_\infty \norma{\v}\norma{\nabla \ph}
    \leq \d \norma{\nabla \ph}^2
    + \cd \norma{\v}^2,
     \\
     \iO (\nabla \ck \star  \ph )\cdot \nabla\ph
     & \leq 
    \norma{\ck}_{W^{1,1}(B_M)}\norma{\ph} \norma{ \nabla \ph}
    \leq \d \norma{\nabla \ph}^2
    + \cd \norma{\ph}^2,
     \\  
	- \iO  (F'' (\ph_1) -F'' (\ph_2) ) \nabla \ph_2 \cdot \nabla \ph 
    & \leq 
    \norma{\ph}_3\norma{\nabla \ph_2}_6\norma{\nabla \ph}
    \\ & 
    \signo{\leq \d \norma{\nabla \ph}^2
    + \cd \norma{\ph}\norma{\nabla \ph}\norma{\ph_2}_{W^{1,6}}^2}
     \\& 
     \leq \d \norma{\nabla \ph}^2
    + \cd \norma{\ph}^2 \norma{\ph_2}_{W^{1,6}}^4.
\end{align*}
In the above computations, we employed integration by parts, the \Holder, Young and Gagliardo--Nirenberg inequalities, as well as the \Lip\ continuity of $F''$ which follows from the separation property \eqref{SP}. 
We then collect the above estimates, selecting $\d$ small enough \last{and}  integrate over time.
Using the regularity \eqref{reg} and applying the Gronwall lemma, we obtain \signo{the first} continuous dependence result \eqref{cont:dep:est1}\signo{.}


\subsubsection*{\signo{Second estimate: part (ii)}}
We now move to the second  continuous dependence estimate. 
\signo{Due to Theorem \ref{ExistCahn}, it follows that}
\begin{align}
  \label{reg1:1}   
  t &\mapsto \norma{\ph_1(t)}^4_{W^{1,6}} +\norma{\ph_2(t)}^4_{W^{1,6}}\in L^1(0,T),
    \quad 
    t \mapsto \norma{\ph_2(t)}^4_{H^2}\in L^1(0,T),
    \\ 
    t & \mapsto \norma{\ph_2(t)}_{W^{1,\infty}}^\an{3} \in L^1(0,T)\signo{.}
    \label{reg1:2}
\end{align}
Then, we take the gradient of \eqref{eq:cd:1} and  the laplacian of \eqref{eq:cd:2} to obtain the corresponding identities
\begin{alignat*}{2}
	& \nabla (\dt\ph)
    + \nabla(\nabla \ph \cdot \v_1)
    {+ \nabla(\nabla \ph_2 \cdot \v)}
    -\nabla \Delta \mu
	=
	0 && \quad \text{in $Q$},
	\\
	& \Delta \mu = 
 -  \div(\nabla\ck \ast \ph) + F''(\ph_1) \Delta \ph 
	+( F''(\ph_1)- F''(\ph_2)) \Delta \ph_2 
    \\ & 
	\qquad\quad
	+ ( (F^{(3)}(\ph_1)- F^{(3)}(\ph_2))|\nabla \ph_2|^2 
    + F^{(3)}(\ph_1)\nabla \ph \cdot (\nabla \ph_1+\nabla\ph_2)	&& \quad \text{in $Q$}.
\end{alignat*}
We test the first one by $\nabla \ph$, the second one by $-\Delta \ph$, and add the resulting equalities leading to a cancellation to infer that
\begin{align*}
	& \frac 12 \frac d{dt} \norma{\nabla \ph}^2
	+ \iO  F''(\ph_1) |\Delta \ph|^2
	=
    \iO \nabla \ph \cdot \v_1 \Delta \ph
    {+ \iO \nabla   \ph_2  \cdot\v\Delta \ph}
	+ \iO \div(\nabla \ck \ast \ph) \Delta \ph
	\\ & \quad 
	- \iO ( F''(\ph_1)- F''(\ph_2)) \Delta \ph_2 \Delta \ph
	- \iO ( (F^{(3)}(\ph_1)- F^{(3)}(\ph_2))|\nabla \ph_2|^2 \Delta \ph \\&\quad
	- \iO F^{(3)}(\ph_1)\nabla \ph \cdot (\nabla\ph_1+\nabla \ph_2)	\Delta \ph
	= \sum_{i=1}^{6}{I}_i.
\end{align*}
Notice that this estimate is indeed rigorous thanks to \rev{\eqref{regul}}. \rev{We now recall the following elliptic inequality
$$
\norma{\ph}_{H^2}\leq C(\norma{\ph}+\norma{\Delta\ph}).
$$}%
\rev{Exploiting it, together with }the \Holder, Young and Agmon inequalities, for a positive \signo{constant} $\d$ yet to be selected, we infer that
\begin{align*}
    I_1 + I_2 
    & \leq 
     {\norma{\nabla \ph}_3\norma{\v_1}_6\norma{\Delta \ph}
    + \norma{\nabla\ph_2}_\infty \norma{\v}\norma{\Delta \ph}}
    \\&\leq\rev{\norma{\nabla\ph}^\frac12\norma{\ph}_{H^2}^\frac12\norma{\v_1}_6\norma{\Delta \ph} + \norma{\nabla\ph_2}_\infty \norma{\v}\norma{\Delta \ph}}
    \\&
    \leq \rev{C\norma{\nabla\ph}^\frac12(\norma{\ph}^\frac12+\norma{\Delta\ph}^\frac12)\norma{\v_1}_6\norma{\Delta \ph} + \norma{\nabla\ph_2}_\infty \norma{\v}\norma{\Delta \ph}}
     \\ & 
   {  \leq 
    {\d \norma{\Delta \ph}^2
    +\cd (\norma{\nabla \ph}^2\norma{\v_1}_6^4
    \rev{+ \norma{\ph}_V^2\norma{\v_1}_6^2}+\norma{ \ph_2}^2_{W^{1,\infty}}\norma{\v}^2),}}
    \\
    I_3 & \leq 
     C \norma{\ph}\norma{ \Delta \ph}
     \leq 
    \d \norma{\Delta \ph}^2
    + \cd \norma{\ph}_V^2,
     \\
    I_4 &  \leq C 
    \norma{\ph}_\infty \norma{\Delta \ph_2}\norma{\Delta \ph}\\&
    \leq \rev{C\norma{\ph}_{V}^\frac12\norma{\ph}_{H^2}^\frac12\Vert \ph_2\Vert_{H^2}\norma{\Delta\ph}}\\&
    \leq \rev{C\norma{\ph}_{V}^\frac12\norma{\ph}^\frac12\Vert \ph_2\Vert_{H^2}\norma{\Delta\ph} + C\norma{\ph}_{V}^\frac12\Vert \ph_2\Vert_{H^2}\norma{\Delta\ph}^\frac32}
    \\&
    \leq \rev{C\norma{\ph}_{V}\Vert \ph_2\Vert_{H^2}\norma{\Delta\ph} + C\norma{\ph}_{V}^\frac12\Vert \ph_2\Vert_{H^2}\norma{\Delta\ph}^\frac32}
    \\ & \leq 
    \rev{\d \norma{\Delta \ph}^2 
    +\cd \norma{\ph}_V^2\norma{\ph_2}_{H^2}^2
    + \cd \norma{\ph}_V^2\norma{\ph_2}_{H^2}^4},
     \\
     I_5 & \leq 
    {C\norma{\ph}_6\norma{
    \nabla \ph_2 }_6^2\norma{ \Delta\ph}}
    \leq 
     \d \norma{\Delta \ph}^2
    +\cd \norma{\ph_2 }_{W^{1,6}}^4\norma{ \ph}^2_V,\\
     I_6 & \leq
    \norma{F^{(3)}(\ph_1)}_\infty(\norma{
    \nabla \ph_1 }_6+\norma{
    \nabla \ph_2}_6)\norma{ \nabla  \ph}_3\norma{ \Delta \ph}
    \\ & 
    \leq 
    \d \norma{\Delta \ph}^2
    + \cd (\norma{\ph_1}^2_{W^{1,6}} + \norma{\ph_2 }_{W^{1,6}}^2)\norma{ \nabla  \ph}\rev{\norma{ \ph}_{{H^2}}}
    \\ &
    \leq 
    {2}\d \norma{\Delta \ph}^2
    + \cd (\rev{\norma{\ph_1}_{W^{1,6}}^2+\norma{\ph_2 }_{W^{1,6}}^2}+\norma{\ph_1}_{W^{1,6}}^4+\norma{\ph_2 }_{W^{1,6}}^4)\norma{   \ph}^2_V . 
\end{align*} 
Observe that,
to handle $I_3$, we have used the inequality
$$
\norma{\div(\nabla\ck\ast \ph)}\leq C\norma{\ph},
$$
which is valid recalling assumption \ref{h4} and \cite[Lemma 2]{Bedrossian}.
Collecting all the estimates above, recalling assumption \ref{h2}, \rev{applying Poincaré's inequality where necessary, recall that we are assuming $(\ph)_\Omega\equiv 0$}, and choosing a consequently small $\d>0$, we end up, after integrating over time for an arbitrary $t \in [0,T)$, with
\begin{align}
    & \nonumber\frac 12 \iO |\nabla \ph(t)|^2
     +\frac{\alpha}{2} \int_{Q_t}|\Delta\ph|^2
     \leq \non
     \frac 12 \iO |\nabla \ph_0|^2
      + C\iot \norma{{\nabla}\ph}^2
     + C \iot \norma{ \ph_2}^2_{W^{1,\infty}}\norma{\v}^2
     \\ & \quad \non
     + \iot C\norma{{\nabla}\ph}^2(\rev{\norma{\ph_2}_{H^2}^2+}\norma{\ph_2}_{H^2}^4)
    \\ & \quad 
    + C \iot (\rev{\norma{\ph_1}_{W^{1,6}}^2+\norma{\ph_2 }_{W^{1,6}}^2}+\norma{\ph_1}_{W^{1,6}}^4+\norma{\ph_2 }_{W^{1,6}}^4{\rev{+\norma{\v_1}^2_6} +\norma{\v_1}^4_6})\norma{ \nabla  \ph}^2.
\label{end}
\end{align}
All the terms are in the usual form in order to apply Gronwall's lemma but the following one that can \an{be} controlled  by \Holder's inequality as
\begin{align*}
     & \iot
    \norma{   \ph_2}_{W^{1,\infty}}^2\norma{\v}^2
    \signo{\leq \iot}\Big(\norma{\ph_2}^2_{W^{1,\infty}}\iO|\v_1-\v_2|^2 \Big)\leq C\norma{\ph_2}^{\an{2}}_{L^\an{3}(0,\signo{t};\Wx{1,\infty})}\norma{\v}_{L^\an{6}(0,\signo{t};\Hs)}^{\an{2}}.
\end{align*}
We then use the regularity \eqref{reg1:1}-\eqref{reg1:2} \last{together with the assumptions on $\v_1,\v_2$}, and apply the Gronwall lemma  to deduce
\eqref{cont:dep:est} concluding the proof.

\end{proof}

\section{The Control Problem}
\label{SEC:CONTROL}
\setcounter{equation}{0}

This section is devoted to the analysis of the optimal control \CP. As mentioned, the control variable consists of a prescribed solenoidal velocity flow occurring in equation \eqref{eq:1}.
From the mathematical properties of the state system \Sys\ addressed in Section \ref{SEC:STATE}, 
the \signo{associated} {\it control-to-state} operator $\S$\signo{,} also referred to as the solution operator, is well-defined and continuous between suitable \signo{Banach} spaces.
For \signo{convenience}, let us repeat here some notation. First, the tracking-type cost functional we aim at minimizing \signo{is defined by}
\eqref{cost}
subject to admissible controls $\v$ \signo{belonging} \an{to}
\begin{align}
    \label{Uad}
    \Vad: = \{ \v \in \L\infty{\LLL^\infty}\cap\L2 {\HHH_\sigma}: \vmin \leq \v \leq \vmax\}.
\end{align}
The specific assumptions on the constants  and target functions in \eqref{cost} and \eqref{Uad} are expressed by \ref{ass:control:4:target}-\ref{ass:control:5:Vad}. 
From now on, as the velocity $\v$ will be constrained in $\Vad$, notice that all the technical requirements in Theorem \ref{ExistCahn} and \ref{THM:CD} are fulfilled.
\last{
Besides, as we will address some  differentiability properties of $\S$, let us take an open ball in the
$L^\infty$-topology containing the set of admissible controls $\Vad$.
Namely, we fix $R>0$ such that 
\begin{align*}
    \Vad\subset \VR := \Big\{ \v \in \L\infty {\LLx\infty} \cap \L2 \Hs: 
    \, \norma{\v}_{\L\infty{\LLL^\infty}} < R \Big\}
    .
\end{align*}
}%
{Then}, the control-to-state operator {$\S$} is the \signo{well-posed} map
\begin{align*}
    \S: {\VR} \to \YY,
    \quad 
    \S : \v \mapsto (\ph,\mu)=(\S_1(\v),\S_2(\v)),
\end{align*}
where $ (\ph,\mu)$ is the unique solution to \Sys\ corresponding to $\v$ and $\YY$ indicates
the {\it state space} which arises from  Theorem \ref{ExistCahn} and \signo{it} is defined as
{
\begin{align}
	\non
	\YY 
	& =
 \YY_1 \times \YY_2
 :=
	\Big(\an{\H1 {\Lx6}}\cap \L\infty V \cap \L4 {\signo{W}} \cap \L3 {\Wx{1,\infty}}\Big)
		\\  \label{def:Y}
		& \qquad 
		\times 
		\Big(\W{1,4} H \cap \C0 V \cap \L8 {\Wx{1,4}} \cap \L4 W\Big).
\end{align}
}%
Besides,  as a consequence of Theorem \ref{THM:CD}, $\S$ is also \Lip\ continuous in the sense expressed \signo{in the} theorem and we also recall that $\ph=\S_1(\v)$ enjoys the separation property \eqref{SP}

The solution operator given above allows us to reduce the optimization problem \CP\ in the usual manner via the {\it reduced cost functional}
\begin{align}
    \label{Jred}
    \Jred(\con) := \JJ (\con;\S_1(\con)),
\end{align}
leading to the minimization problem
\begin{align*}
    \min_{\con \in \Vad} \Jred(\con).
\end{align*}

\subsection{Existence of optimal controls}
The first step consists in proving Theorem \ref{THM:EXCONTROL}. This can be done by a straightforward application of the direct method of calculus of variations.
\begin{proof}[Proof of Theorem \ref{THM:EXCONTROL}]
   First, we notice that $\J$ is bounded from below as it is nonnegative.
   Let $\{\v_n \}_n \subset \Vad$ be a minimizing sequence for $\Jred$ and let $(\ph_n,\mu_n ):=\SS(\v_n)$ denote the sequence of the corresponding states, $n\in\mathbb{N}$. Then, up to nonrelabeled subsequences, there exists $\vopt\in \Vad$ such that
   \begin{align}
\v_n\to \vopt\quad \text{weakly* in }\signo{\L\infty{\LLL^\infty} \cap \L2 {\Hs}}.
\label{weakconv}
\end{align}
By the results of Theorem \ref{ExistCahn}, since the sequence $\{(\ph_n,\mu_n )\}$ is associated to the same initial datum $\ph_0$ and $\v_n$ is uniformly \signo{bounded}, we immediately infer the following uniform bounds
{
   \begin{align*}
    &
    \norma{\ph_n}_{\YY_1 \cap L^\infty(Q)} 
    + 
    \Vert \mu_n\Vert_{\YY_2}
    \leq C.
   \end{align*}
}%
This implies, by standard weak and weak$^*$ compactness arguments, that there exists $(\bph,\bmu)$ such that, up to subsequences, 
{
\begin{align*}
    &\ph_n\to \bph\quad\text{weakly* in }\YY_1\cap L^\infty(Q),
    \quad 
    \mu_n\to \bmu\quad\text{weakly* in }\YY_2,\\&
\end{align*}}%
which imply, by the Aubin--Lions--Simon lemma, that
{
\begin{align*}
    \varphi_n\to \bph\quad\text{strongly in } C^0([0,T];H^r(\Omega)) \cap L^4(0,T;W^{s,2}(\Omega)),  \quad \forall s\in [0,2),r\in [0,1).
\end{align*}}%
Note now that, for any $n\in\mathbb{N}$, $(\ph_n,\mu_n)$ satisfies
\eqref{weak-nCH-2:1}-\eqref{weak-nCH-2:2} with $\v:=\v_n$. From these convergences and \eqref{weakconv} we can easily pass to the limit as $n\to\infty$ in the weak formulation \eqref{weak-nCH-2:1}-\eqref{weak-nCH-2:2} and deduce that $(\bph,\bmu)$ is a weak solution according to the definition of Theorem \ref{ExistCahn} with $\v:=\vopt$. Therefore, by the uniqueness of weak solutions, we immediately infer $\SS(\vopt)=(\bph,\bmu)$, i.e., the pair $(\vopt,\bph)$ is admissible for \signo{the minimization problem} $\CP$. By weak lower sequential semicontinuity of \signo{norms}, it readily follows that $(\vopt,\bph)$ is optimal for $\JJ$. Indeed, recall that the above convergences also imply $\ph_n(T)\to \bph(T)$ strongly in $H$ as $n \to \infty$. Therefore, $(\vopt,\bph)$ yields a solution to $\CP$ and the proof is concluded.  
\end{proof}
\subsection{Differentiability properties of the solution operator}
The natural subsequent step is providing some optimality conditions for the minimizers of \CP. 
Since the set of admissible controls $\Vad$ is convex, it is well-known that the first-order optimality conditions are characterized by
a suitable variational inequality of the form 
\begin{align}\label{abs:varineq}
    \<D\Jred (\opt), \v- \opt> \geq 0
    \quad \forall \v \in \Vad,
\end{align}
where $\Jred$ is the reduced cost functional introduced above, and $D\Jred$ stands for its \Fre\ derivative in a suitable mathematical framework.
The aim of this section is to set the ground to rigorously justify the above formula.
In this direction, the first step is to obtain some differentiability properties of the control-to-state operator $\S$.


We then fix a control $\vopt \in {\VR}$ 
and denote by $\optstate\signo{:=\S(\opt)}$ the corresponding state.
Then, the linearized system to \Sys\ at $\opt$, for any $\w \in \L2 {\Hs}$, reads (in \signo{strong} form) as 
\begin{alignat}{2}
	\label{eq:lin:1}
	& \dt \xi + \nabla \xi \cdot \vopt + \nabla \bph \cdot \w - \Delta \eta
	=
	0
	\quad && \text{in $Q$,}
	\\
	\label{eq:lin:2}
	& \eta = - \ck \ast \xi + F'' (\bph)\xi
	\qquad && \text{in $Q$,}
	\\
	\label{eq:lin:3}
	& \dn \eta = 
	0
	\quad && \text{on $\Sigma$,} 
	\\
	\label{eq:lin:4}
	& \xi(0)=0
	\quad && \text{in $\Omega$.}
\end{alignat}
\Accorpa\Lin {eq:lin:1} {eq:lin:4}
The weak well-posedness of the \signo{above} system follows.
\begin{theorem}
    \label{THM:LIN}
    Suppose that \ref{h2}-\ref{h3} and \ref{ass:control:1:kernel}-\ref{ass:control:2:initialdata} and \ref{ass:control:5:Vad} are fulfilled. Then, for every $\w \in \L2 \Hs$, there exists a unique solution $(\xi,\th)$ to the linearized system \Lin\ \signo{in the sense that}  
    \begin{align*}
        \xi & \in \H1 \Vp \cap \C0 H \cap \L2 V ,
        \quad 
        \th  \in \L2 V,
    \end{align*}
    and it \signo{fulfills}
    \begin{align}
        \non
        & \<\dt \xi, v >
        - \iO  \xi\, \vopt \cdot \nabla v 
        + \iO  \nabla \bph \cdot \w v
        + \iO \nabla  \eta \cdot \nabla v
        =0
        \\ \label{wf:lin:1}
        & \qquad \text{for every $v \in V $, and \aet},
        \\ \label{wf:lin:2}
        & 
        \iO \eta v 
        = - \iO \ck \ast \xi v 
        + \iO F''(\bph)\xi v
         \quad \text{for every $v \in V $, and \aet},
    \end{align}
    as well as the initial condition
    \begin{align*}
        \xi(0)=0 \quad \text{in $\Omega$.}
    \end{align*}
\end{theorem}

\begin{proof}
The proof of existence is standard and can be performed by approximation argument, for instance, using a Faedo--Galerkin  scheme.
For this reason, we proceed formally by avoiding the introduction of any approximation
scheme and limiting ourselves to  provide formal estimates.

\step 
First estimate

We test \eqref{eq:lin:1} by $\NN \xi$ (observe that $ \xi_\Omega=0$), \eqref{eq:lin:2} by $-\xi$ and add the resulting identities to obtain 
\begin{align*}
    \frac 12  \frac d{dt} \norma{\xi}^2_*
    + \iO F''(\bph) |\xi|^2
    = 
   \iO  \xi \,\vopt \cdot \nabla \NN \xi
    - \iO  (\nabla \bph \cdot \w) \NN \xi
    + \iO (\ck \ast \xi )\xi.
\end{align*}
From \ref{h2}\signo{, we readily} get
\begin{align*}    
    \iO F''(\bph) |\xi|^2 \geq \alpha\norma{\xi}^2.
\end{align*}
Let us  now move to bound the terms on the \rhs. In this direction, we use the definition of $\NN$ introduced in Section \ref{SUBSEC:NOT}, the \Holder\ and  Young inequalities, and the regularity of $\bph$ as solution to \Sys\ in the sense of Theorem \ref{ExistCahn}.
Recalling that $\bph\in[-1,1]$, we have that
\begin{align*}
    \iO  \xi \vopt \cdot \nabla \NN \xi & \leq \norma{\xi} \norma{\opt}_\infty \norma{\xi}_*
    \leq\frac{\alpha}6\norma{\xi}^2 + C \norma{\xi}^2_*,
    \\
     - \iO  (\nabla \bph \cdot \w) \NN \xi &=\iO\bph(\w\cdot\nabla\last{\NN}\xi) \leq \norma{ \bph}_\infty\norma{\w}\norma{\nabla\NN\xi}
     \leq C\norma{\w}^2+C \norma{\xi}^2_*,
     \\ 
     \iO (\ck \ast \xi )\xi & = \<\ck \ast \xi ,\xi> \leq \norma{\ck \ast \xi}_V\norma{\xi}_*
     \leq 
     C \norma{\xi} \norma{\xi}_*\leq\frac{\alpha}6 \norma{\xi}^2+ C \norma{\xi}^2_*,
\end{align*}
where we also owe to 
$$
\norma{\ck\ast \xi}_V\leq C(\norma{\ck\ast \xi}+\norma{\nabla\ck\ast \xi})\leq C\norma{\xi},
$$
which follows from Young's inequality for convolutions and \ref{ass:control:1:kernel}.
Integrating over time and employing Gronwall's lemma, recalling the regularity on $\w$, yield
\begin{align}
    \norma{\xi}_{\L\infty \Vp \cap \L2 H}
    \signo{\leq C \norma{\w}_{\L2 {\Hs}}}
    \leq C. 
    \label{above_estimate}
\end{align}

\step 
Second estimate

Using the above estimate \eqref{above_estimate}, it is not difficult to realize, from a comparison argument in \eqref{eq:lin:2}, that also
\begin{align*}
   \norma{\eta}_{\L2 H}
    \signo{\leq C \norma{\w}_{\L2 {\Hs}}}
    \leq C.    
\end{align*}

\step 
Third estimate 

Next, we test \eqref{eq:lin:1} by $\xi$.
Then, we consider the gradient of equation \eqref{eq:lin:2} and test it by $-\nabla \xi$. Adding the identities leads us to 
\begin{align*}
    & \frac 12  \frac d{dt} \norma{\xi}^2
    + \iO F''(\bph)|\nabla \xi|^2
    =
    - \iO  \nabla \xi \cdot \vopt \xi
        - \iO  \nabla \bph \cdot \w \xi
    \\ & \quad     \rev{+} \iO (\nabla \ck \star \xi) \cdot \nabla \xi
        - \iO  F^{(3)}(\bph) \nabla \bph \, \xi \cdot \nabla \xi.
\end{align*}
The second term on the left-hand side is bounded from below\signo{, again,} due to \ref{h2}. Next, we observe that, being \last{$\vopt\in \VR$} divergence-free,
$$
 - \iO  \nabla \xi \cdot \vopt \xi
 =- \iO  \nabla  \Big(\frac {\xi^2} 2\Big)  \cdot \vopt 
 =  
 \iO   \Big(\frac {\xi^2} 2\Big)  \div \vopt 
 - \iG \frac {\xi^2} 2 \,  \vopt \cdot  \nnn
 =
 0.
$$
To bound the other terms, we owe to the \Holder,\ Gagliardo--Nirenberg and Young inequalities.  Namely,  we find that
\begin{align*}
     \rev{\iO} (\nabla \ck \star \xi) \cdot \nabla \xi
    & \leq \norma{\ck}_{W^{1,1}(B_M)} \norma{\xi}\norma{\nabla \xi} 
     \leq \frac {\alpha}8 \norma{\nabla \xi}^2
     + C \norma{\xi}^2,
     \\
     - \iO  \nabla \bph \cdot \w \xi
   &  = \iO \bph \w   \cdot \nabla \xi \leq \norma{\bph}_\infty \norma{\w}\norma{\nabla \xi}
    \leq \frac {\alpha}8 \norma{\nabla \xi}^2
     + C\norma{\w}^2 ,
     \\
     - \iO  F^{(3)}(\bph) \nabla \bph \,\xi \cdot \nabla \xi & \leq 
     \norma{F^{(3)}(\bph)}_\infty
     \norma{ \nabla \bph}_6
     \norma{ \xi}_3
     \norma{\nabla  \xi}
    \\ & 
    \leq \frac {\alpha}8 \norma{\nabla \xi}^2
     +C \norma{  \bph}_{W^{1,6}}^2
     \norma{ \xi}\norma{\nabla \xi}
     \\ & 
    \leq \frac {\alpha}4 \norma{\nabla \xi}^2
     +C \norma{  \bph}_{W^{1,6}}^4
     \norma{ \xi}^2.
\end{align*}
We then integrate over time, recalling also that 
 $t \mapsto \norma{  \bph(t)}_{W^{1,6}}^4 \in L^1(0,T)$ and $\w \in \L2 \Hs$, and  apply the Gronwall's lemma to conclude that 
\begin{align*}
    \norma{\xi}_{\L\infty H \cap \L2 V}
    \signo{\leq C \norma{\w}_{\L2 {\Hs}}}
    \leq C.    
\end{align*}
As \signo{before}, comparison in \eqref{eq:lin:2} readily  produces
\begin{align*}
   \norma{\eta}_{\L2 V}
   \signo{\leq C \norma{\w}_{\L2 {\Hs}}}
    \leq C.    
\end{align*}

\step 
Fourth estimate

From a comparison argument in \eqref{eq:lin:1} it is a standard matter to infer that
\begin{align*}
    \norma{\dt \xi}_{\L2 \Vp}
    \signo{\leq C \norma{\w}_{\L2 {\Hs}}}
    \leq C.    
\end{align*}
Finally, the continuity property of $\xi$ can be easily deduced due to the continuous embedding $\H1  \Vp \cap \L2 V \emb \C0 H.$

As far as uniqueness is concerned, we notice that \signo{system} \Lin\ is linear.  Thus, the above estimates \signo{readily entail} uniqueness.
Indeed, for two special solutions $(\xi_i,\eta_i)$, $i=1,2$, we set $(\xi,\eta):=(\xi_1-\xi_2,\eta_1-\eta_2)$ and observe that the above estimates hold for $(\xi,\eta)$ with $C=0$\signo{. H}ence, the claimed uniqueness follows.
\end{proof}

It is now naturally expected that, provided we select the correct Banach spaces, the linearized
system captures the behavior of the \Fre\ derivative of the solution operator, that is, the identity $D\S(\opt)\signo{[\w]} = (\xi, \eta)$ \signo{holds} in a suitable mathematical setting.
This is rigorously stated in the following result, where the following Banach space appears
\begin{align*}
    {\cal X}:= \big(\signo{\H1 \Vp \cap \C0 H} \cap \L2 V \big) \times
    \L2 V.
\end{align*}

\begin{theorem}
    \label{THM:FRE}
   Suppose that \ref{h2}-\ref{h3} and \ref{ass:control:1:kernel}-\ref{ass:control:2:initialdata} and \ref{ass:control:5:Vad} are in force. Then, the solution mapping $\S$ is \Fre\ differentiable at every $\opt \in \VR$ as a mapping from ${\L6 \Hs}$ into $\cal X$.
    In addition, it holds that 
    \begin{align*}
        \text{$D\S(\opt) \in {\cal L} ({\L6 \Hs}, {\cal X})$
        \,\, 
        and
        \,\,
        $D\S(\opt)[\w] =(\xi,\eta)$}
    \end{align*}    
     with $(\xi,\eta)$ being the unique solution to the linearized system \Lin\ corresponding to $\w$ as given by Theorem \ref{THM:LIN}.
\end{theorem}
\begin{proof}[Proof of Theorem \ref{THM:FRE}]
Suppose \signo{the identity} $D\S(\opt)[\w] =(\xi,\eta)$ \signo{has been proven}. Then, it readily follows from Theorem \ref{THM:LIN} that $D\S(\opt) \in {\cal L} ({\L6 \Hs}, {\cal X})$ \signo{as $\w\mapsto (\xi,\eta)$ is linear and continuous from $\L2 \Hs $ to $\cal X$}.

\noindent Let us then show that $D\S(\opt)[\w] =(\xi,\eta)$ by checking that 
\begin{align}
     \label{fre:formal:1}
    \frac{\norma{\S(\opt + \w) - \S(\opt) - (\xi, \eta)}_{\cal X}}{\norma{\w}_{{\L6 \Hs}}} \to 0 
    \quad \text{as $\norma{\w}_{{\L6 \Hs}}\to 0.$ }
\end{align}
As a matter of fact, what we are going to prove is a specific estimate that will imply the above property. Namely, we aim at showing that  
\begin{align*}
    {\norma{\S(\opt + \w) - \S(\opt) - (\xi, \eta)}_{\cal X}} \leq  C {\norma{\w}_{{\L6 \Hs}}^2}
\end{align*}
for a suitable positive constant $ C>0$. 
Without loss of generality, we tacitly assume from now on that the norm of the increment \signo{$\w$} is small enough so that $\opt + \w$ remains in the open set $\VR$.
Upon setting
\begin{align*}
    (\hph, \hmu):= \S( \vopt + \w), 
    \quad 
    \optstate:=\S(\opt),
    \quad 
    \psi := \hph - \bph - \xi, 
    \quad 
    \th := \hmu - \bmu - \eta,
\end{align*}
the above inequality \eqref{fre:formal:1} amounts proving the \an{existence} of a constant $C>0$ such that
\begin{align*}
    \norma{(\psi, \th)}_{\cal X} \leq C {\norma{\w}^{2}_{{\L6 \Hs}}}.
\end{align*}
In the direction of checking the above estimate, we write the system solved by $(\psi, \th)$, noticing that, by the previous regularity results, we have
\signo{$(\psi,\th) \in{ \cal X}$.}
Namely, we consider the difference between \Sys\ considered at $\opt +\w$ and  at $\opt$, together with \Lin. We thus infer that the pair $(\psi, \th)$ is a (weak) solution to 
\begin{alignat}{2}
	\label{eq:fre:1}
	& \dt \psi + \nabla \psi \cdot \vopt + \nabla (\hph - \bph) \cdot \w - \Delta \th
	=
	0
	\quad && \text{in $Q$,}
	\\
	\label{eq:fre:2}
	& \th = - \ck \ast \psi + [F' (\hph)-F' (\bph)-F'' (\bph)\xi]
	\qquad && \text{in $Q$,}
	\\
	\label{eq:fre:3}
	& \dn \th = 
	0
	\quad && \text{on $\Sigma$,} 
	\\
	\label{eq:fre:4}
	& \psi(0)=0
	\quad && \text{in $\Omega$.}
\end{alignat}
\Accorpa\Fresys {eq:fre:1} {eq:fre:4}
Before proceeding, it is worth noticing that, as a consequence of 
Theorem \ref{ExistCahn}\signo{,} recall $\opt \in \VR$ and \last{\ref{ass:control:1:kernel}- \ref{ass:control:2:initialdata} are in force}\signo{,} we have that
\begin{align*}
    & \exists \, \beta\in (0,1): \quad 
    \ph \in C^{\beta,\frac \beta2}(\ov Q),
    \quad 
    \text{and}
    \quad 
   \last{ \norma{\ph}_{ L^8(0,T;W^{1,4}(\Omega)) \cap \L4 {\Hx2} \cap \L3 {\Wx{1,\infty}}}}
   \leq C,
\end{align*}
with $\ph\in\{\hph,\bph\}$, as well as the validity of the strict separation property 
\begin{align*}
    \exists \, \delta>0 : \quad \max\{|\hph(x,t)|,|\bph(x,t)|\} < 1 - \delta \quad  \forall (x,t)\in \last{\ov{Q}}.
\end{align*}
Thus, using also \ref{ass:control:6:potreg}, there exists a positive constant $C$ such that
\begin{align*}
    \big|{F^{(k)}(\hph(x,t))-F^{(k)}(\bph(x,t))}\big|
    \leq C |{\hph(x,t) - \bph(x,t)}|
    \quad 
    \forall (x,t)\in \last{\ov{Q}},
    \quad k=0,...,3.
\end{align*}
Besides, from Theorem \ref{THM:CD}, it holds the continuous dependence result
\begin{align} \label{cd:fre}
    \norma{\hph-\bph}_{\L\infty V \cap \L2 W} \leq C\norma{\w}_{{\L6 \Hs}}.
\end{align}
Finally, we recall Taylor's formula with integral remainder for a generic regular function $f\in C^2(\erre)$:
\begin{align}
	\label{taylor:f}
	f(\hph)-f(\bph)-f'(\bph)\xi \,=\, f'(\bph) \psi + \RR (\hph-\bph)^2,
\end{align}
where the remainder $\RR$ is defined as
\begin{align*}
	\RR:= \int_0^1 f''(\bph+s (\hph-\bph)) (1-s)\, {\rm ds}.
\end{align*}
For our purposes, $f = F^{(k)}$ with $k=1,2$ and we abuse notation by denoting the corresponding remainders with the same symbol $\RR$. Due to the regularity assumed in \ref{h2}, along with the separation property \eqref{sepa}, it holds that the associated remainders are uniformly bounded.
We are now ready to show the claimed estimate.

\step 
First estimate

Below, we are going to consider the gradient of equation \eqref{eq:fre:2}.
Hence, we highlight the identity
\begin{align*}
    & \nabla [F' (\hph)-F' (\bph)-F'' (\bph)\xi]
     = F'' (\hph) \nabla \hph-F'' (\bph) \nabla \bph-F^{(3)} (\bph)\nabla \bph \, \xi - F'' (\bph)\nabla \xi
    \\ 
    & \quad  =
    [F'' (\hph) -F'' (\bph) -F^{(3)} (\bph)\xi] \,\nabla \bph
    \\ & \qquad 
    + (F''(\hph) - F''(\bph))\nabla (\hph - \bph)
    + F''(\bph)\nabla \psi
    =: {\Theta}+ F''(\bph)\nabla \psi.
\end{align*}
Next, we test \eqref{eq:fre:1} by ${\cal N} \psi + \psi $ (observe that $\psi_\Omega=0$), \eqref{eq:fre:2} by $\psi$, the gradient of \eqref{eq:fre:2} by $\nabla\psi$, and add the resulting equalities to obtain that
\begin{align*}
    & \frac 12 \frac d {dt} (\norma{\psi}^2_*+\norma{\psi}^2)
    + \iO F''(\bph) |\nabla \psi|^2
    \\ & \quad 
    = 
    - \iO \nabla \psi \cdot \opt \NN \psi 
    - \iO \nabla (\hph - \bph) \cdot \w  (\NN \psi + \psi)
    + \iO (\ck \ast \psi) \psi 
     \\ & \qquad 
    + \iO [F' (\hph)-F' (\bph)-F'' (\bph)\xi] \psi
    + \iO (\nabla \ck \star \psi)\cdot \nabla  \psi 
    - \iO \Theta \cdot \nabla \psi
   = \sum_{i=1}^{6} I_i,
\end{align*}
where we recall that, arguing as previously done,
$$
    \iO (\nabla\psi\cdot \vopt)\psi=0.
$$
Using similar computations as in the other proofs, for a positive \signo{constant} $\d$ yet to be chosen, we infer that
\begin{align*}
    I_1 +I_3 {+ I_5}& \leq 
    \norma{\psi}\norma{\opt}_\infty \norma{\nabla\NN\psi}
    + {C}\norma{ \ck}_{W^{1,1}(B_M)}\norma{\psi}( \norma{\psi} + \norma{\nabla \psi})
    \\ & 
    \leq \d \norma{\nabla \psi}^2
    + \cd (\norma{\psi}^2_*+\norma{\psi}^2),
    \\
    I_2 & \leq 
    \norma{\nabla (\hph- \bph)}_3\norma{\w}(\norma{\NN\psi}_6+\norma{\psi}_6)
    \leq \d \norma{ \nabla\psi}^2
    +  {\cd}\norma{\nabla (\hph- \bph)}_3^2\norma{\w}^2
     \\ & \quad 
     {+ \cd (\norma{\psi}_*^2+\norma{\psi}^2)}
    \\ &  \leq \d \norma{ \nabla\psi}^2
    +  {\cd}\norma{\nabla (\hph- \bph)}(\rev{\norma{\hph-\bph}+}\norma{\Delta (\hph- \bph)})\norma{\w}^2
    {+ \cd (\norma{\psi}_*^2+\norma{\psi})^2}
    \\ & 
    \leq \d \norma{ \nabla\psi}^2
    + \norma{\w}^4 
    + {\cd}\norma{\hph- \bph}^{2}_{\L\infty V}(\rev{\norma{ \hph- \bph}^2_V+}\norma{\Delta (\hph- \bph)}^2) 
    \\&
    \quad {+ \cd (\norma{\psi}_*^2+\norma{\psi}^2)},
    \\ I_4 & \leq
    \norma{F''(\bph)}_\infty\norma{\psi}^2
    + \norma{\RR}_\infty \norma{\hph-\bph}_4^2\norma{\psi}
    \leq C (\norma{\psi}^2 + \norma{\hph-\bph}^4_V),
\end{align*}
where we applied also Poincaré's inequality and Taylor's formula \eqref{taylor:f} with $f = F'$ and integrated by parts the first integral  $I_1$.
As for $I_6$, using the above definition of $\Theta$, we infer that 
\begin{align*}
    I_6  & \leq 
    \norma{\Theta}\norma{\nabla \psi}
    \leq \d \norma{\nabla \psi}^2
    + \cd \norma{\Theta}^2
    \\ & \leq 
    \d \norma{\nabla \psi}^2
    + \cd \norma{F^{(3)}(\bph)}_\infty ^2\norma{\psi}^2_3\norma{\nabla \bph}^2_6
    + \cd \norma{\RR}_\infty^2\norma{\hph-\bph}^4_6\norma{\nabla \bph}^2_6
    \\ & \leq 
    \d \norma{\nabla \psi}^2
    + \cd \norma{\psi}^2_3\norma{ \bph}^2_{W^{1,6}}
    + \cd \norma{\hph-\bph}^4_V\norma{ \bph}^2_{W^{1,6}}
    \\ & \leq 
    2\d \norma{\nabla \psi}^2
    + \cd \norma{\psi}^2\norma{ \bph}^4_{W^{1,6}}
    + \cd \norma{\hph-\bph}^4_V\norma{ \bph}^2_{W^{1,6}},
\end{align*}
owing to Taylor's formula \eqref{taylor:f} with $f = F''$.
Collecting all the above estimates, recalling \ref{h2} and choosing $\varepsilon$ suitably small, we get
\begin{align*}
     \frac 12 \frac d {dt} (\norma{\psi}^2_*+\norma{\psi}^2)
    + \alpha\norma{\nabla \psi}^2
    &\leq C(\norma{\psi}^2_*+\norma{\psi}^2)+ \norma{\hph-\bph}^4_V
    + \norma{\w}^4 
    \\& \quad 
    + {C}\norma{\hph- \bph}^{2}_{\L\infty V}(\rev{\norma{ \hph- \bph}^2_V+}\norma{\Delta (\hph- \bph)}^2).
\end{align*}
We now apply Gronwall's lemma, noticing that 
$t \mapsto \norma{ \bph(t)}^4_{W^{1,6}} \in L^1(0,T)$ and recalling the stability estimate \eqref{cd:fre} and obtain that
\begin{align*}
    \norma{\psi}_{\L\infty H \cap \L2 V}
    \leq C \norma{\w}^2_{{\L6 \Hs}}.
\end{align*}

\step 
Second estimate

From the above estimate it readily follows, from a comparison argument in \eqref{eq:fre:2}, that 
\begin{align*}
    \norma{\th}_{\L2 V}
    \leq C \norma{\w}^2_{{\L6 \Hs}}.
\end{align*}

\step
\signo{Third estimate}

\signo{
We finally go back to equation \eqref{eq:lin:1} to infer by comparison, using the above bounds, that
\begin{align*}
    \norma{\dt \psi}_{\L2 \Vp}
    \leq C \norma{\w}^2_{{\L6 \Hs}}.
\end{align*}
This latter also entails $  \norma{\psi}_{\C0 H} \leq C \norma{\w}^2_{{\L6 \Hs}}$ due to standard compactness embeddings.
}%
This concludes the proof as the above estimates imply \eqref{fre:formal:1}, whence the claim.
\end{proof}

\subsection{Optimality conditions}
The final step consists in proving Theorem \ref{THM:FOC}, but first we point out the following intermediate result.

\begin{theorem}
\label{THM:FIRST:FOC}
    Assume that \ref{ass:control:1:kernel}-\ref{ass:control:6:potreg} are in force. Let $\vopt\in \Vad$ be an optimal control with corresponding state $\optstate=\S(\opt)$.
    Then, it holds that 
    \begin{align}
	\label{foc:first}
    \alphaQ \intQ (\bph - \ph_Q)\xi 
    + \alphaO \iO (\bph(T) -\ph_\Omega)\xi(T)
    + \alphav\intQ \opt \cdot (\v-\vopt)
	\geq 0
	\quad \forall \con \in \Vad,
    \end{align}
where $\lin$ stands for the unique linearized variables associated to $\w = \v - \opt$ as given by Theorem \ref{THM:LIN}.
\end{theorem}

\begin{proof}[Proof of Theorem \ref{THM:FIRST:FOC}]
    This is a straightforward consequence of the abstract result \eqref{abs:varineq} along with Theorem \ref{THM:LIN}, the definition of the cost functional $\J$, and the chain rule.
\end{proof}

As customary, \eqref{foc:first} \signo{is} not helpful in numeric schemes and have to be somehow simplified.  This is done by the help of an additional problem related  to \Sys\ called adjoint system. It consists of a backward-in-time parabolic system in the variables $\adj$ and, in its strong form, it reads as
\begin{alignat}{2}
	\label{eq:adj:1}
	& -\dt p  - \ck \ast q + F'' (\bph)q - \nabla p\cdot \vopt 
	=
	\alphaQ (\bph-\ph_Q)
	\quad && \text{in $Q$,}
	\\
	\label{eq:adj:2}
	& q= - \Delta p
	\qquad && \text{in $Q$,}
	\\
	\label{eq:adj:3}
	& \dn p = 
	0
	\quad && \text{on $\Sigma$,} 
	\\
	\label{eq:adj:4}
	& p(T)=\alphaO(\bph(T)-\ph_\Omega)
	\quad && \text{in $\Omega$.}
\end{alignat}
\Accorpa\Adj {eq:adj:1} {eq:adj:4}

\begin{theorem}
    \label{THM:ADJ}
   Suppose that \ref{ass:control:1:kernel}-\ref{ass:control:6:potreg} hold.
   Let $\opt$ be an optimal control with corresponding state $\optstate\signo{:=\S(\opt)}$. 
   Then, there exists a unique solution $\adj$ to the adjoint system \Adj\ in the sense that  
    \begin{align*}
        p & \in \H1 \Vp \cap \C0 {{V}} \cap \L2 {{W}} ,
        \quad
        q \in \L2 H,
    \end{align*}
    and it solves
    \begin{align}
        \non
        & -\<\dt p, v >
        - \iO (\ck \ast q) v
        +\iO  F'' (\bph)q v
        - \iO (\nabla p\cdot \vopt )v
	=
	\iO \alphaQ (\bph-\ph_Q) v
        \\ \label{wf:adj:1}
        & \qquad 
        \text{for every $v \in V $, and \aet},
        \\ \label{wf:adj:2}
        & 
        \iO q v 
        = \iO\nabla p \cdot \nabla v
         \quad \text{for every $v \in V $, and \aet},
    \end{align}
    as well as the terminal condition
    \begin{align*}
        p(T)=\alphaO(\bph(T)-\ph_\Omega) \quad \text{in $\Omega$.}
    \end{align*}
\end{theorem}
\begin{proof}[Proof of Theorem \ref{THM:ADJ}]
Again, we proceed formally. The rigorous proof can be performed, e.g., by using \signo{a} Faedo--Galerkin approximation scheme.

\step
First estimate

We test \eqref{eq:adj:1} by $q$, \eqref{eq:adj:2} by $\dt p$ and add the resulting identities to infer that
\begin{align*}
    - {\frac 12 \frac d{dt} \norma{\nabla p}^2
   } +  \iO F''(\bph) |q|^2
    =
     \iO (\ck \ast q) q
     + \iO (\nabla p\cdot \vopt )q
     +\iO \alphaQ (\bph-\ph_Q)q.
\end{align*}
As done several times, we owe to the convexity of $F$ in \ref{h2} to derive that the \signo{second} term on the \lhs\ is bounded from below by $\alpha \norma{q}^2$.
As for the terms on the \rhs, we use basic computations \an{by means of Young's inequality for convolutions} to obtain
\begin{align*}
     \iO (\ck \ast q) q
      =
     \<\ck \ast q, q>&\an{\leq \norma{\ck\ast q}_V\norma{q}_*}
     \leq 
     \frac{\alpha}6 \norma{q}^2
     + C \norma{q}^2_*
     \leq 
     \frac{\alpha}6 \norma{q}^2
     + C \norma{\nabla p}^2,
     \\
     \iO (\nabla p\cdot \vopt )q
     +\iO \alphaQ (\bph-\ph_Q)q
     & \leq 
     \norma{\nabla p} \norma{\opt}_\infty\norma{q}
     +\an{C}(\norma{\bph}_\infty +\norma{\ph_Q})\norma{q}
     \\ & \leq 
     \frac{\alpha}3 \norma{q}^2
     + C (1+\norma{\nabla p}^2+\norma{\ph_Q}^2) .
\end{align*}
\last{
Observe that, exploiting \eqref{wf:adj:2}, in the first estimate we have used the \signo{bound}
$$
\norma{q}_*\leq \norma{\nabla p}.
$$
}
Next, we integrate over $(t,T)$ for an arbitrary $t \in [0,T)$ to infer that 
\begin{align*}
    \norma{\nabla p(t)}^2
    + \alpha \norma{q}^2
    \leq 
    \norma{\nabla p(T)}^2
    + C \int_t^T |\nabla p|^2
    + C.
\end{align*}
Then, we recall \eqref{eq:adj:4}, \ref{ass:control:4:target} and the fact that {$\bph \in \C0 V$} as strong solution in the sense of Theorem \ref{ExistCahn}, whence the terminal condition on the \rhs\ of the above inequality can be bounded.
The backward-in-time Gronwall lemma then entails that 
\begin{align*}
    \norma{p}_{\L\infty V}
    + \norma{q}_{\L2 H}
    \leq C.
\end{align*}

\step
Second estimate

A comparison argument in \eqref{eq:adj:2} and elliptic regularity theory produce
\begin{align*}
    \norma{p}_{\L2 W}
    \leq C.
\end{align*}

\step
Third estimate

Finally, it is a standard matter to derive from \eqref{eq:adj:1} that
\begin{align*}
    \norma{\dt p}_{\L2 \Vp} \leq C.
\end{align*}

\last{In conclusion}, we recall that \Adj\ is linear: thus, the above computations already entail uniqueness.
\end{proof}
\begin{proof}[Proof of Theorem \ref{THM:FOC}]
To prove the theorem, we compare the variational inequalities \eqref{foc:final} with \eqref{foc:first}.
We then realize that it suffices to check that 
\begin{align}
    \label{foc:identity}
    - \intQ \PPP_\sigma(p \nabla \bph) \cdot (\v - \opt)
    = \alphaQ \intQ (\bph - \ph_Q)\xi
    + \alphaO \iO (\bph(T) - \ph_\Omega)\xi(T)
\end{align}
with $\xi$ being the first component of the unique solution to \Lin\ associated to $\w= \v - \opt$.
In this direction, we multiply \eqref{wf:lin:1} by $p$ with $\w=\v- \opt$, \eqref{wf:lin:2} by $-q$ and integrate over $(0,T)$.
We then owe to the well-known integration-by-parts formula
for functions belonging to
$\H1 \Vp \cap \L2 H$
as well as to \eqref{wf:adj:1} and \eqref{wf:adj:2}
to obtain that    
\begin{align*}
        0& = \iot \<\dt \xi,p> 
        - \intQ \xi \opt \cdot \nabla p
        + \intQ  \nabla \bph \cdot (\v- \opt) p
        + \intQ \nabla \xi \cdot \nabla p
       \\ & \quad 
       + \intQ [-\eta  - \ck \ast \xi + F'' (\bph)\xi]q
       \\ &  =
       -\iot \<\dt p,\xi> 
       + \intQ \xi[  - \ck \ast q + F'' (\bph)q - \nabla p\cdot \vopt ]
       + \intQ \eta[- q - \Delta p]
        \\ & \quad 
        +\iO p(T) \xi(T)
        -\iO p(0) \xi(0)
       +\intQ p\nabla \bph \cdot( \v- \opt).
    \end{align*}
Using the initial and terminal conditions in
\eqref{eq:lin:4} and \eqref{eq:adj:4}, we are led to \eqref{foc:identity}.
In particular, we notice that, since $\signo{\w=}\v-\vopt \in \Hs$ and the projector $\PPP_\sigma$ is selfadjoint, it holds that
$$
    \intQ p\nabla \bph \cdot( \v- \opt)
    =
    \intQ \PPP_\sigma (p\nabla \bph) \cdot( \v- \opt),
$$
concluding the proof.
\end{proof}


\section*{Acknowledgments}
\rev{We thank the anonymous referees for their valuable comments and remarks, which significantly improved the clarity of our work.}
Partial support 
from the MIUR-PRIN Grant 
2020F3NCPX ``Mathematics for industry 4.0 (Math4I4)'' is  gratefully acknowledged.
Besides, the authors are affiliated to the GNAMPA (Gruppo Nazionale per l'Analisi Matematica, 
la Probabilit\`a e le loro Applicazioni) of INdAM (Isti\-tuto 
Nazionale di Alta Matematica).
\rev{They also acknowledge that the
present research has been supported by MUR, grant Dipartimento di Eccellenza 2023-2027.}


\vspace{3truemm}

\Begin{thebibliography}{10} \footnotesize

	
\bibitem{Adams} Adams, R., Fournier, J., \textit{Real Interpolation of Sobolev Spaces on Subdomains of $\mathbb{R}^n$}. Canadian Journal of Mathematics, \textbf{30} (1978), 190-214.

\bibitem{num5} Allen, S., Cahn, J.W., \textit{A microscopic theory for antiphase boundary motion and its application to antiphase domain coarsening}, Acta Metall., \textbf{27} (1979), 1084-1095.

\bibitem{Bedrossian}  Bedrossian, J., Rodrìguez, N., Bertozzi, A., \textit{Local and global well-posedness for an
aggregation equation and Patlak--Keller--Segel models with degenerate diffusion}, Nonlinearity, \textbf{24} (2011), 1683-1714.

\bibitem{BH}
Brangwynne, C. P., Hyman, A. A.,  
\textit{ Beyond
stereospecificity: liquids and mesoscale
organization of cytoplasm}, Dev. Cell, \textbf{21} (2011), 14-16.

\bibitem{BM}
Brezis, H., Mironescu, P.,
{\it Where Sobolev interacts with Gagliardo–Nirenberg}, 
J. Funct. Anal., {\bf 277} (2019), 2839-2864.


\bibitem{num13} Cahn, J.W., Hilliard, J.E., \textit{Free energy of a non-uniform system. I. Interfacial free energy},
J. Chem. Phys., \textbf{28} (1958), 258-267.

\rev{
\bibitem{CES}
Carrillo, J.A., Elbar, C.,  Skrzeczkowski, J.,
{\it Degenerate Cahn--Hilliard systems: from nonlocal to local}, (2023), Preprint arXiv:2303.11929.
}

\bibitem{CM} Cherfils, L., Miranville, A.,  Zelik, S., \textit{The Cahn--Hilliard equation with logarithmic potentials}, Milan J. Math., \textbf{79} (2011), 561-596.

{
    \bibitem{CGS1}
    Colli, P., Gilardi, G., Sprekels, J.,
    {\it Optimal velocity control of a viscous Cahn--Hilliard system with convection and dynamic boundary conditions},
    SIAM J. Control Optim., {\bf 56} (2018), 1665-1691.
}

{
    \bibitem{CGS2}
    Colli, P., Gilardi, G., Sprekels, J.,
    {\it Optimal velocity control of a convective Cahn--Hilliard system with double obstacles and dynamic boundary conditions: a `deep quench' approach},
    J. Convex Anal., {\bf 26} (2019), 485-514.
}

{
\bibitem{CSS4}
    Colli, P., Signori, A., Sprekels, J.,
    {\it Optimal control problems with sparsity for phase field tumor growth models involving variational inequalities},
    J. Optim. Theory Appl., (2022). doi.org/10.1007/s10957-022-02000-7.
}

{
    \bibitem{CL}
    Cristini, V., Lowengrub, J.,
    {\it Multiscale Modeling of Cancer: An Integrated Experimental and Mathematical},
    Modeling Approach. Cambridge University Press, Leiden (2010).
}

\bibitem{Scarpa1}  Davoli, E., Scarpa, L., Trussardi, L., \textit{Nonlocal-to-local convergence of Cahn--Hilliard equations: Neumann boundary conditions and viscosity terms}, Arch. Ration. Mech. Anal. \textbf{239} (2021), 117-149.

\bibitem{Scarpa2} Davoli, E., Scarpa, L., Trussardi, L., \textit{Local asymptotics for nonlocal convective Cahn--Hilliard equations with $W^{1,1}$ kernel
and singular potential}, J. Differential Equations, \textbf{289} (2021), 35-58.


\bibitem{Dolgin}  Dolgin, E., \textit{What lava lamps and vinaigrette can teach us about cell biology}, Nature, \textbf{555} (2018), 300-302.

\rev{
\bibitem{ES}
Elbar, C.,  Skrzeczkowski, J.,
{\it Degenerate Cahn--Hilliard equation: from nonlocal to local},
Journal of Differential Equations,
{\bf  364} (2023), 576-611.
}

\bibitem{EG} Elliott, C.M., Garcke, H., \textit{On the Cahn--Hilliard equation with degenerate mobility}, SIAM J. Math. Anal., \textbf{27} (1996), 404-423.

\rev{
\bibitem{F}
Fornoni, M.,
{\it Optimal distributed control for a viscous non-local tumour growth model}, 
App. Math. Opt., {\bf 89}, Online first (2024), doi:10.1007/s00245-023-10076-4.
}


{
    \bibitem{Fmob}
    Frigeri, S., 
    {\it On a nonlocal Cahn--Hilliard/Navier--Stokes system with degenerate mobility and singular potential for incompressible fluids with different densities},
    Ann. Henri Poincaré,
    {\bf 38} (2021), 647-687.
}

\bibitem{Frig} Frigeri, S., Gal, C.G., Grasselli, M.,
	\textit{Regularity results for the nonlocal Cahn--Hilliard equation with singular potential and degenerate mobility} (2021),
	Journal of Differential Equations, \textbf{287},  295-328.
 
 \bibitem{Frigeri} Frigeri, S., Gal, C.G., Grasselli, M., \textit{On nonlocal Cahn--Hilliard--Navier--Stokes systems in two dimensions}, J. Nonlinear Science, \textbf{26} (2016), 847-893.

{
    \bibitem{FGGS}
    Frigeri, S., Gal, C.G., Grasselli, M., Sprekels, J.,
    {\it Two-dimensional nonlocal Cahn--Hilliard--Navier--Stokes systems with variable viscosity, degenerate mobility and singular potential}, Nonlinearity, {\bf 32},  (2019), 678.
}

\bibitem{FG} Frigeri, S, Grasselli, M., \textit{Nonlocal Cahn--Hilliard--Navier--Stokes systems with singular potentials}, Dyn. Partial Diff. Eqns., \textbf{9} (2012), 273-304.

\bibitem{FG1} Frigeri, S, Grasselli, M., \textit{Global and trajectory attractors for a nonlocal Cahn--Hilliard--Navier--Stokes system}, J. Dynam. Differential Equations, \textbf{24} (2012), 827-856.

{
    \bibitem{FGS}
     Frigeri, S., Grasselli, M., Sprekels, J.,
    {\it Optimal distributed control of two-dimensional nonlocal Cahn--Hilliard--Navier--Stokes systems with degenerate mobility and singular potential},
    App. Math. Opt., {\bf 81(3)} (2020) 899-931.
}

 \bibitem{FGK}
 Frigeri, S., Grasselli, M., Krej\acc{c}\'i, P.,
 {\it Strong solutions for two-dimensional nonlocal Cahn--Hilliard--
Navier--Stokes systems},
{ J. Differential Equations}, {\bf 255} (2013), 2587-2614.

{
    \bibitem{FLS}
    Frigeri, S., Lam, K.F., Signori, A.,
    {\it Strong well-posedness and inverse identification problem of a non-local phase field tumor model with degenerate mobilities. }
    European J. Appl. Math., {\bf 33}(2) (2022), 267-308.
}

{
    \bibitem{FRS}
   Frigeri, S., Rocca, E., Sprekels, J.,
    {\it Optimal distributed control of a nonlocal Cahn--Hilliard/Navier--Stokes system in two dimensions},
    SIAM Journal on Control and Optimization, {\bf 54(1)} (2016), 221-250.
}

\bibitem{GGG0} 
{ Gal, C.G.,  Giorgini A., Grasselli, M.}, 
{\textit{ The nonlocal Cahn--Hilliard equation with singular potential: well-posedness, regularity and strict separation property}}, 
{J. Differential Equations},
 {\textbf{263}} (2017), {5253-5297}.

\bibitem{GGG} Gal, C.G., Giorgini, A., Grasselli, M., \textit{The separation property for 2D Cahn--Hilliard equations: local, nonlocal and fractional energy cases}, Discrete Contin. Dyn. Syst., \textbf{43} (2023),
2270-2304.

\bibitem{AGGP} Gal, C.G., Giorgini, A., Grasselli, M., Poiatti, A., \textit{Global well-posedness and convergence
to equilibrium for the Abels-Garcke-Grün model with nonlocal free energy} (2023), \rev{J. Math. Pures Appl.
(9) \textbf{178}, 46-109.}

\rev{
\bibitem{GPf} 
Gal, C.G., Poiatti, A., {\it Unified framework for the separation property in binary phase
segregation processes with singular entropy
densities}, Researchgate preprint  (2023), 10.13140/RG.2.2.35972.30089/1.
}

\bibitem{GLS_OPT}
    Garcke, H., Lam,  K.F., Signori, A.,
    {\it Sparse optimal control of a phase field tumour model with mechanical effects},
    SIAM J. Control Optim., {\bf 59}(2) (2021), 1555-1580.

\bibitem{giacomin} Giacomin, G., Lebowitz, J.L., \textit{Exact macroscopic description of phase segregation in model alloys with long range
	interactions}, Phys. Rev. Lett., \textbf{76} (1996), 1094-1097.

\bibitem{giacomin1} Giacomin, G., Lebowitz, J.L., \textit{Phase segregation dynamics in particle systems with
	long range interations}, I. Macroscopic limits, J. Stat. Phys.,  \textbf{87} (1997), 37-61.

\bibitem{giacomin2} Giacomin, G., Lebowitz, J.L., \textit{Phase segregation dynamics in particle systems with
	long range interations}, II. Interface motion, SIAM J. Appl. Math., \textbf{58} (1998), 1707-1729.	

\bibitem{GGM}  Giorgini, A., Grasselli, M., Miranville, A., T\textit{he Cahn--Hilliard--Oono equation with singular potential}, Math. Models Methods
Appl. Sci., \textbf{27} (2017), 2485-2510.

\bibitem{St} Giorgini, A., \textit{On the separation property and the global attractor for the nonlocal Cahn--Hilliard equation in three dimensions} (2023), Preprint arXiv:2303.06013.

{
   \bibitem{HW}
   Hinterm\"uller, M., Wegner, D., 
   {\it Optimal control of a semidiscrete Cahn--Hilliard--Navier--Stokes system}, SIAM Journal on Control and Optimization, {\bf 52} (2014), 747-772.
}

\bibitem{Lady}  Lady\acc{z}enskaja, O.A., Solonnikov, V.A., Ural’ceva, N.N., \textit{Linear and quasilinear equations of parabolic type}, AMS Transl. Monographs 23, AMS, Providence, R.I. 1968.



\bibitem{M} Miranville, A., \textit{The Cahn--Hilliard Equation: Recent Advances and Applications},
CBMS-NSF Regional Conf. Ser. in Appl. Math., SIAM, Philadelphia, PA., 2019.


\bibitem{Muramatu} Muramatu, T., \textit{On Besov Spaces of Functions Defined in General Regions}, Publications of The Research Institute for Mathematical Sciences, \textbf{6} (1970), 515-543.


\bibitem{P} Poiatti, A., \textit{The 3D strict separation property for the nonlocal Cahn--Hilliard
equation with singular potential} (2022), \rev{Anal. PDE, to appear (see also: Preprint arXiv:2303.07745).}

\bibitem{Pruss} Pr\"uss, J., Schnaubelt, R.,
\textit{Solvability and Maximal Regularity of Parabolic Evolution Equations with Coefficients Continuous in Time},
Journal of Mathematical Analysis and Applications, \textbf{256} (2001), 405-430.

{
    \bibitem{RSS}
     Rocca, E., Scarpa, L., Signori, A.,
    {\it Parameter identification for nonlocal phase field models for tumor growth via optimal control and asymptotic analysis},
    Math. Models Methods Appl. Sci., {\bf 31}(13) (2021), 2643-2694.
}

{
    \bibitem{RS}
    Rocca, E., Sprekels, J., 
    {\it Optimal distributed control of a nonlocal convective Cahn--Hilliard equation by the velocity in three dimensions},
    SIAM Journal on Control and Optimization,
    {\bf 53(3)} (2015) 1654-1680.
}

{
\bibitem{SS}
    Scarpa, L., Signori, A.,
     {\it On a class of non-local phase-field models for tumor growth with possibly singular potentials, chemotaxis, and active transport},
     Nonlinearity,  {\bf 34} (2021), 3199-3250. 
}

{
    \bibitem{ST}
    Sprekels, J., Tröltzsch, F.,
    {\it Sparse optimal control of a phase field system with singular potentials arising in the modeling of tumor growth},
    ESAIM Control Optim. Calc. Var., {\bf 27},
     (2021) S26.
}

{
    \bibitem{ST2}
    Sprekels, J., Tröltzsch, F.,
    {\it Second-order sufficient conditions for sparse optimal control of singular Allen--Cahn systems with dynamic boundary conditions},
    \rev{Discrete Contin. Dyn. Syst. S, Online first (2023),  doi: 10.3934/dcdss.2023163.}
}

\bibitem{Triebel} Triebel, H., \textit{Interpolation Theory, Function Spaces, Differential Operators}, 2. Rev. and enl. ed. Heidelberg: Barth, 1995.

\bibitem{Wu} Wu, H., \textit{A review on the Cahn–Hilliard equation: classical results and recent advances in dynamic boundary conditions}, J. Electronic Research Archive, \textbf{30} (2022), 2788-2832.

{
    \bibitem{ZL}
    Zhao, X., Liu, C.,
    {\it Optimal control for the convective Cahn--Hilliard equation in 2D case},  App. Math. Opt., {\bf 70} (2014), 61-82.
}

\End{thebibliography}

\End{document}

\bye